\documentclass[numbers,webpdf,imaiai]{ima-authoring-template}%

\graphicspath{{Fig/}}

\theoremstyle{thmstyletwo}%

\usepackage{soul}

\newtheorem{theorem}{Theorem}[section]
\newtheorem{proposition}[theorem]{Proposition}%
\newtheorem{assumption}[theorem]{Assumption}
\newtheorem{lemma}[theorem]{Lemma}
\newtheorem{corollary}[theorem]{Corollary}
\newtheorem{claim}[theorem]{Claim}

\newtheorem{remark}[theorem]{Remark}%

\newtheorem{definition}[theorem]{Definition}

\numberwithin{equation}{section}

\begin{document}

\DOI{DOI HERE}
\copyrightyear{2022}
\vol{00}
\pubyear{2022}
\access{Advance Access Publication Date: Day Month Year}
\appnotes{Paper}
\copyrightstatement{Published by Oxford University Press on behalf of the Institute of Mathematics and its Applications. All rights reserved.}
\firstpage{1}


\title[Concentration Inequalities for Sums of Markov Dependent Random Matrices]{Concentration Inequalities for Sums of Markov Dependent Random Matrices}

\author{Joe Neeman, Bobby Shi*, and Rachel Ward
\address{\orgname{The University of Texas at Austin} }}

\authormark{Joe Neeman, Bobby Shi*, and Rachel Ward}

\corresp[*]{Corresponding author: \href{bhshi@utexas.edu}{bhshi@utexas.edu}}

\received{Date}{0}{Year}
\revised{Date}{0}{Year}
\accepted{Date}{0}{Year}


\abstract{We give Hoeffding and Bernstein-type concentration inequalities for the largest eigenvalue of sums of random matrices arising from a Markov chain.  We consider time-dependent matrix-valued functions on a general state space, generalizing previous results that had only considered Hoeffding-type inequalities, and only for time-independent functions on a finite state space.  In particular, we study a kind of noncommutative moment generating function, provide tight bounds on this object, and use a method of Garg et al. to turn this into tail bounds.  Our proof proceeds spectrally, bounding the norm of a certain perturbed operator.  In the process we make an interesting connection to dynamical systems and Banach space theory to prove a crucial result on the limiting behavior of our moment generating function that may be of independent interest.}

\maketitle

\section{Introduction}
The study of concentration inequalities has now become textbook material, with a variety of applications \cite{boucheron2013}.  Two of the most widely used and studied inequalities for scalar-valued random variables in the independent setting are Hoeffding's inequality \cite{hoeffding1963}, which operates under a boundedness assumption, and Bernstein's inequality \cite{bernstein1927}, which operates under an additional bounded variance assumption.

With the wide usage of these inequalities, various generalizations have been made.  One line of work has sought to relax the independence assumption, deriving concentration inequalities for sums of Markov-dependent random variables.  Starting with the result of Gillman \cite{gillman1993}, improvements and refinements have been made to address the general setting of time-independent functions of a nonreversible Markov chain on a general state space \cite{lezaud1998, lezaud2001, paulin2015, leon2004, fan2021, jiang2018, healy2008}, resulting in Hoeffding and Bernstein-type concentration inequalities that generalize nicely the independent setting; these results are largely spectral in nature.  

Another line of work has sought to generalize the scalar setting to concentration inequalities of sums of independent random matrices, beginning with the work of Ahlswede and Winter \cite{ahlswede2002}.  Several authors have followed the Ahlswede-Winter approach to develop analogs of Hoeffding's and Bernstein's inequalities for sums of random matrices; some works include \cite{gross2011, recht2011, Christofides2008, oliveira2009, oliveira2010}.  However, results using this framework often have variance proxies that are suboptimal.  To this end, \cite{tropp2011, tropp2015} developed techniques to circumvent this barrier, allowing for concentration inequalities that are much tighter in many cases.  These results are powerful and easy to use \cite{tropp2016, TROPP2018}, with various applications \cite{avner2011, DRINEAS2011, lopez2014}.

We combine the above two lines of study in our work to develop concentration inequalities -- Hoeffding and Bernstein-type -- for the largest eigenvalue of sums of Markov dependent Hermitian random matrices, by combining the Ahlswede-Winter style argument for random matrices with the spectral techniques used to study Markov dependent scalar-valued random variables.  The starting point is a deep result of Garg et al. \cite{garg2018}, which provides a multi-matrix Golden-Thompson inequality that supports a spectral approach.  Our results are broad and general: we provide inequalities for \textit{time-dependent} matrix-valued functions of a \textit{nonreversible} Markov chain on \textit{general (continuous) state spaces}.  In this way, our Hoeffding-type bounds vastly generalize previous work and provide sharper constants, and our Bernstein-type bounds are new in the literature (Section \ref{sec:final}).  Along the way we prove a novel Perron-Frobenius type limit lemma for the matrix setting (Lemmas \ref{lem:limit}, \ref{lem:simple_function}),  necessitated by our assumptions of nonversibility of the chain on continuous state spaces, that may be of independent interest.

Thus, this work can be seen as a first systematic study of concentration inequalities in the Markov-dependent matrix setting.  Each of our results involve two bounds: one bounding a noncommutative moment generating function of the form
\[\mathbb{E}\bbb{\norm{\prod_{j=1}^n \exp\bb{\frac{\theta}{2}X_j}}_F^2}\]
and one bounding the tail probabilities
\[\operatorname{Pr}\bb{\lambda_{\text{max}}\bb{\sum_{j=1}^n X_j}\geq t}\]
for Hermitian random matrices $X_1, \dots, X_n$ arising as functions of a Markov chain.  We expect our results on the former to be optimal, directly generalizing the scalar setting.  

These results are readily applicable to a wide variety of settings.  For example, the best known guarantees for offline principal component analysis (PCA) are derived from a combination of the standard matrix Bernstein's inequality along with Wedin's perturbation theorem \cite{wedin1972, jain2018, kumar2024}.  Now, if the sample matrices are not independent, but instead arise from a Markov chain, we again have a bound for leading eigenvector estimation using our bounds.  Our results are also directly applicable to sums of matrix-valued random variables sampled from a random walk on an expander; Theorem \ref{thm:markov_matrix_hoeffding} improves on the best known results in this setting \cite{wigderson2005randomness, garg2018}.  There is also a line of work in machine learning that realizes stochastic gradient descent (SGD) with constant step size as a Markov chain with a stationary distribution \cite{bach2013, Dieuleveut2020, mandt2017}.  Additionally, Hessian information has been used widely to give improved algorithms and better understand model performance \cite{ghorbani2019, li2020}.  Our results are able to give precise bounds on the convergence of the largest eigenvalue of the empirical mean of the Hessian matrices along the path of SGD to the largest eigenvalue of the expected Hessian with respect to the stationary measure.  This has many interesting and relevant applications in statistical learning and generalization.  More broadly, our bounds are useful in the study of online algorithms, where matrix-valued data is received and processed in a streaming fashion from an underlying Markov chain.  In Section \ref{sec:pca} we more closely focus on an application of our main results to offline PCA for samples arising as a function of an underlying Markov chain; we present a bound that is a direct generalization of the best known bound for offline PCA in the i.i.d. setting.  This result was first stated in \cite{kumar2024}, where to our knowledge it was the first in the literature to give this sort of bound for offline Markov PCA ; it directly uses our Theorem \ref{thm:markov_matrix_bernstein}.

\subsection{Related work}

In terms of results, the literature on concentration inequalities for sums of Markov-dependent random matrices is fairly small, in contrast to the scalar setting.  Recent progress in this area builds on the work of \cite{sutter2017multivariate, garg2018}, which develop a powerful Golden-Thompson inequality \cite{golden1965lower, thompson1965inequality} that allows for a spectral analysis of a useful moment generating function; the latter provides a corresponding Hoeffding's inequality for when the stationary distribution is uniform and the chain is stationary.  The work of \cite{qiu2020} extends this to arbitrary stationary distributions and initial distributions; these works are valid only for finite state spaces and reversible chains.  Our work generalizes these results to general state spaces and nonreversible chains using significantly different techniques, which allows us to also improve the constants in the bounds.  Furthermore, our Hoeffding-type result is in terms of both the largest and smallest eigenvalues, matching the type given in the scalar setting; and we give a Bernstein's inequality, which is new in the literature.  

In terms of techniques, most related to ours are the works of {\cite{leon2004, MIASOJEDOW2014, fan2021, jiang2018, lezaud1998, paulin2015, healy2008}}, which give Hoeffding's and Bernstein's inequalities for sums of Markov-dependent \textit{scalar} random variables.  Broadly speaking, the idea is to bound the largest norm of a perturbed operator.  We prove the corresponding spectral properties of a Markov operator that is the tensored form of the operator appearing in these works in Section \ref{sec:operator}; in particular, we generalize to continuous state spaces using related functional analytic ideas.  However, previous techniques are not sufficient for us to prove our main results; our proof of the result on the limit of the moment generating function (Lemma \ref{lem:limit}) takes a detour through dynamical systems theory and Banach space theory.

In contrast to the Markov-dependent setting, the sums of independent random matrices have been studied in depth, and is the model for which we give our bounds.  In particular, our results mirror the type given in \cite{ahlswede2002, Christofides2008, gross2011} in terms of variance proxy, and the analysis proceeds from the same starting point: bounding a matrix moment generating function; in contrast, later works of \cite{oliveira2010, tropp2011} sharpen the variance proxy.  Our analysis of the moment generating function arising from the matrix-valued Golden-Thompson inequality result in bounds resembling more of the former; we believe our analysis is likely sharp given the form of the moment generating function, so an improvement in terms of the variance proxy would likely require a different approach and is an interesting open problem.

Besides the spectral approach, there have been a myriad of ways to derive concentration inequalities relaxing the independence assumption, both in the scalar and matrix settings.  The papers of \cite{Glynn2002hoeffding, Adamczak2015exponential} exploit regeneration-type minorization conditions to derive exponential concentration for ergodic sums; though these works do not assume a nonzero spectral gap they often have less explicit constants.  The work of \cite{paulin2015} uses Marton couplings to obtain concentration inequalities for a scalar-valued function of a single random sample.  The papers of \cite{paulin2016, mackey2014} leverage Efron-Stein inequalities and exchangeability to develop concentration inequalities for random matrices, also relaxing the independence assumption.  Similarly, recent works of \cite{aoun2020matrix, kathuria2020matrix,huang2021} demonstrate that matrix functional inequalities, i.e., Poincar\'e, directly translate to matrix concentration inequalities;\, \cite{kathuria2020matrix} develops a Bernstein's inequality for strongly Rayleigh distributions;\,\cite{kyng2018} gives a Chernoff bound for a similar setting.  These largely address single sample concentration of (Lipschitz) functionals; however, some methods can be generalized to product measures \cite{huang2021nonlinear}, recovering some of the Efron-Stein results.  And the work of \cite{garg2018} sketches a method to reduce studying concentration of random variables sampled from a Markov chain to concentration of sums of martingale random variables \cite{chung2006concentration, tropp2011}; though the resulting bounds are suboptimal, they show that qualitatively concentration for Markov chains is more generic than just the scalar or matrix settings (see \cite{ledoux1991series}).  These techniques may be very useful in improving the variance proxy, better matching the independent setting, where a direct spectral analysis may prove insufficient; in particular, functional inequalities have emerged as a powerful general tool to derive concentration.  It may be interesting to see if this can be demonstrated in the scalar setting as well.

\section{Main Results}\label{sec:results}
In this section we give our main results.  Each will hold under a combination of the following assumptions.

\begin{assumption}\label{ass:markov_chain}
$P$ is a discrete-time Markov chain on a continuous state space $\mathcal{X}$ with stationary distribution $\mmu$ and absolute spectral gap $\lambda$ (Definition \ref{def:absolute_spectral}).  $s_1, \dots, s_n$ is a sequence of states driven by $P$ with initial distribution $\mmu$.
\end{assumption}

\begin{assumption}\label{ass:real_F}
    $F_j:\mathcal{X}\to\bbR^{d\times d}$, $j=1, \dots, n$, is a sequence of functions each mapping to real symmetric $d\times d$ matrices.  Each $F_j$ is $\ell_2(\mmu\otimes \mathbf{1})$ measurable (see Section \ref{sec:preliminaries}).
\end{assumption}

\begin{assumption}\label{ass:hoeffding}
    For all $j$, $\mathbb{E}_{\mmu}[F_j(x)]=0$ and $a_jI\preceq F_j(x)\preceq b_jI$ for all $x\in \mathcal{X}$.
\end{assumption}

\begin{assumption}\label{ass:bernstein}
    For all $j$, $\mathbb{E}_{\mmu}[F_j(x)]=0$ and there exists a constant $\mathcal{V}_j$ such that $\norm{\mathbb{E}_{\mmu}\bbb{F_j(x)^2}}\leq \mathcal{V}_j$ for all $x\in \mathcal{X}$.  Moreover, there exists an absolute constant $\mathcal{M}$ such that $\norm{F_j(x)}\leq \mathcal{M}$ for all $j, x\in \mathcal{X}$.
\end{assumption}

The following two theorems each consist of two statements.  The first bounds a certain type of moment generating function arising from the expected Frobenius norm of a product of matrix exponentials.  The second statement turns this into a tail bound.  We only give upper tail bounds; the corresponding statement for lower tail bounds are clear and the proof proceeds in almost the exact same way.

\begin{theorem}[Markov Matrix Hoeffding, Real]\label{thm:markov_matrix_hoeffding}
Instate Assumptions \ref{ass:markov_chain}, \ref{ass:real_F}, \ref{ass:hoeffding}.  Then for any $\theta>0$, 
\[
    \mathbb{E}_{\mmu}\bbb{\norm{\prod_{j=1}^n \exp\bb{\frac{\theta}{2}F_j(s_j)}}_F^2}\leq d\exp\bb{\frac{\theta^2}{2}\cdot \alpha(\lambda)\cdot \frac{\sum_{j=1}^n (b_j-a_j)^2}{4}}
\]
and for any $t>0$, 
\[
    \operatorname{Pr}_{\mmu}\bb{\lambda_{\text{max}}\bb{\sum_{j=1}^n F_j(s_j)}\geq t}
        \leq d^{2-\pi/4}\exp\bb{\frac{-t^2/(8/\pi^2)}{\alpha(\lambda)\cdot \sum_{j=1}^n (b_j-a_j)^2}}.
\]
Here $\alpha(\lambda)=(1+\lambda)/(1-\lambda)$, where $\lambda$ is the absolute spectral gap (Definition \ref{def:absolute_spectral}).
\end{theorem}
The first statement in the above result reveals the sub-Gaussian nature of the moment generating function; in the scalar case, it exactly reduces to the standard moment generating function.  Using our methods, this bound matches the scalar case given in \cite{fan2021}.  The second statement provides a large deviation tail bound; the bound on the lower eigenvalue follows analogously.  Our results offer four main improvements to that of \cite{qiu2020}: first, we give strictly better constants, both in terms of $d$ and the absolute constants in the exponent (the $\alpha(\lambda)$ differs slightly from their equivalent term for the mixing of the Markov chain, but is strictly better and is more classical); second, we allow for time-dependent functions $F_j$, which is a strict generalization; third, our inequality is in terms of both a lower and upper eigenvalue bound (the $a_j, b_j$ parameters), giving a result that is directly comparable to the classic Hoeffding's lemma; and fourth, our bounds apply to general state spaces.

\begin{theorem}[Markov Matrix Bernstein, Real]\label{thm:markov_matrix_bernstein}
    Instate Assumptions \ref{ass:markov_chain}, \ref{ass:real_F}, \ref{ass:bernstein}.  Let $\sigma^2=\sum_{j=1}^n \mathcal{V}_j$.  Then for any $0<\theta<\log(1-\lambda)/\mathcal{M}$,
    \[\mathbb{E}_{\mmu}\bbb{\norm{\prod_{j=1}^n \exp\bb{\frac{\theta}{2}F_j(s_j)}}_F^2}\leq d\exp\bb{\frac{\sigma^2}{\mathcal{M}^2}\bb{e^{\mathcal{M}\theta}-\mathcal{M}\theta-1}+\frac{\sigma^2}{\mathcal{M}^2}\cdot \frac{\lambda (e^{\mathcal{M}\theta}-1)^2}{1-\lambda e^{\mathcal{M}\theta}}}.\]
    Furthermore, if $\theta<(1-\lambda)/(8\mathcal{M}/\pi)$, then for any $t> 0$,
    \[\operatorname{Pr}_{\mmu}\bb{\lambda_{\text{max}}\bb{\sum_{j=1}^n F_j(s_j)}\geq t}\leq d^{2-\pi/4}\exp\bb{\frac{-t^2/(32/\pi^2)}{\alpha(\lambda)\cdot \sigma^2+\beta(\lambda)\cdot \mathcal{M}t}}.\]
    Here
    \[\alpha(\lambda)=\frac{1+\lambda}{1-\lambda},\quad \beta(\lambda)=\begin{cases}
    \frac{4}{3\pi},\quad\,\, \lambda=0,\\
    \frac{8/\pi}{1-\lambda},\quad 0<\lambda<1
    \end{cases},\]
    where $\lambda$ is the absolute spectral gap (Definition \ref{def:absolute_spectral}).
\end{theorem}
This is to our knowledge the first result giving a Bernstein-type inequality for sums of random matrices arising from markov chains.  The first statement in the above result is made up of two terms; the $(e^{\mathcal{M}\theta}-\mathcal{M}\theta-1)\sigma^2/\mathcal{M}^2$ term coincides with the convex function that makes up the classic Bernstein's and Bennett's inequality.  The second term in the exponent reflects the influence of the Markov chain; compared to previous scalar results \cite{jiang2018, paulin2015}, our techniques are able to recover a slightly improved version of this second term in terms of absolute constants.  A linear algebraic perspective of the first statement above is that the first term bounds how much the product of matrix exponentials increases the magnitude of vectors in the direction of $\mathbf{1}$, which is the leading eigenvector of the operator $P$.  The second term is a bound on how much the product of matrix exponentials increases the magnitude of orthogonal directions, hence involving $\lambda$.  

We now give two corollaries of the above results, generalizing the assumptions above; the proofs are in the appendix and are classical.  Our first corollary generalizes to complex Hermitian matrices via a complexification technique \cite{dongarra1984eigenvalue}:
\begin{corollary}[Extension to Complex Matrices]\label{cor:complex_matrix_concentration}
    Instate Assumption \ref{ass:markov_chain}.  Now assume that $F_j:\mathcal{X}\to\bbC^{d\times d}$, $j=1, \dots, n$, is a sequence of functions each mapping to complex Hermitian $d\times d$ matrices. 

    Under Assumption \ref{ass:hoeffding} (resp. Assumption \ref{ass:bernstein}), the conclusion of Theorem \ref{thm:markov_matrix_hoeffding} (resp. Theorem \ref{thm:markov_matrix_bernstein}) holds with an extra multiplicative factor of $2$ on the right-hand side.
\end{corollary}

Our second corollary generalizes to the case that the Markov chain starts at a distribution $\boldsymbol{\nu}$:
\begin{corollary}[Extension to Nonstationary Chains]\label{cor:nonstationary}
    Let $P$ be a Markov chain on general state space $\mathcal{X}$ with stationary distribution $\mmu$ and absolute spectral gap $\lambda$.  Let $s_1, \dots, s_n$ be a sequence of states driven by $P$ with initial distribution $\boldsymbol\nu$, where $\boldsymbol\nu \ll \mmu$.  Instate Assumption \ref{ass:real_F}.  Let $\operatorname{ess\,sup}\frac{d\boldsymbol{\nu}}{d\mmu}$ be the essential supremum of the Radon-Nikodym derivative $\frac{d\boldsymbol{\nu}}{d\mmu}$.

    Under Assumption \ref{ass:hoeffding} (resp. Assumption \ref{ass:bernstein}), the conclusion of Theorem \ref{thm:markov_matrix_hoeffding} (resp. Theorem \ref{thm:markov_matrix_bernstein}) holds with an extra multiplicative factor of $\operatorname{ess\,sup}\frac{d\boldsymbol{\nu}}{d\mmu}$ on the right-hand side.
\end{corollary}
\begin{remark}
    In many applications these bounds are used to determine how many samples are needed (i.e., how long the chain has to run) in order for the tail probability to be less than some fixed value; in these cases the number of samples will depend only logarithmically on $\operatorname{ess\,sup}\frac{d\boldsymbol{\nu}}{d\mmu}$.  Alternatively, the probability of an event under the Markov chain driven by $P$ initialized at $\boldsymbol{\nu}$ can be bounded by the probability of an event under the Markov chain driven by $P$ initialized at $\mmu$ up to an extra multiplicative factor of $1+\chi^2\infdivx{\boldsymbol{\nu}}{\mmu}$; see Lemma 29 of \cite{kook2024covarianceestimationusingmarkov}.
\end{remark}

We now give some remarks on the above theorems.  First, for the first statements regarding the moment generating function in both Theorems \ref{thm:markov_matrix_hoeffding} and \ref{thm:markov_matrix_bernstein}, the variance proxies are asymptotically optimal for a class of Markov chains.  Let $\mmu$ be a distribution on $\mathcal{X}$ and let $P$ be a transition kernel such that $P(x, y)=\lambda \mathbbm{1}_{x=y}+(1-\lambda)\mmu(y)$, so that $\mmu$ is the stationary distribution.  Let $f:\mathcal{X}\to\mathbb{R}$ be any scalar-valued function, and let $F:\mathcal{X}\to\mathbb{R}^{d\times d}$ be defined as $F(x)=f(x)I$, so that these are all real symmetric.  It is straightforward that
\begin{equation}\label{eq:variance_remark_mgf}
    \mathbb{E}_{\mmu}\bbb{\norm{\prod_{j=1}^n \exp\bb{\frac{\theta}{2}F(s_j)}}_F^2}=d\mathbb{E}_{\mmu}\bbb{\exp\bb{\theta \sum_{j=1}^n f(s_j)}}
\end{equation}
since all $F(x)$ commute.  

Under this setup, start with our Bernstein assumptions.  Let $f$ such that $\mathbb{E}_{\mmu}[f]=0$, $\mathbb{E}_{\mmu}[f^2]= \sigma^2$, and $|f|\leq \mathcal{M}$.  A central limit theorem of \cite{geyer1992} states that for a two state Markov chain driven by a $P$ of our above form that
\begin{equation}\label{eq:var_asy}
\mathcal{V}_{\text{asy}}:=\lim_{n\to\infty} \operatorname{Var}\bb{\frac{1}{\sqrt{n}}\sum_{j=1}^n f(s_j)}=\frac{1+\lambda}{1-\lambda}\cdot \sigma^2,
\end{equation}
which is called the asymptotic variance.  Now suppose that $\mathcal{V}$ is such that $\mathbb{E}_{\mmu}[f^2]\leq \mathcal{V}$; it is classical that any such variance proxy $\mathcal{V}$ that satisfies $g(\theta)=n\mathcal{V}\theta^2/2+o(\theta^2)$ for a function $g$ that upper bounds the scalar moment generating function via 
\[\mathbb{E}_{\mmu}\bbb{\exp\bb{\theta \sum_{j=1}^n f(s_j)}}\leq e^{g(\theta)}\]
upper bounds $\operatorname{Var}\bb{\frac{1}{\sqrt{n}}\sum_{j=1}^n f(s_j)}$ \cite{jiang2018}.  This lower bound is indeed obtained by our Markov kernel $P$ and choice of function $f$, and therefore the same holds for our matrix-valued function $F$.

Now assume $f$ are Rademacher, i.e., the probability under $\mmu$ that $f(x)=1$ is $1/2$ and $f(x)=-1$ is $1/2$.  Define $F$ the same way as above, and so we see that $\mathbb{E}_{\mmu}[F]=0$ and $-I\preceq F(x)\preceq I$ for all $x\in \mathcal{X}$.  Via Equation \ref{eq:variance_remark_mgf} we again reduce to the scalar case since $F$ commute, and through \ref{eq:var_asy} we have that $\mathcal{V}_{\text{asy}}=(1+\lambda)/(1-\lambda)$, since the variance of $f$ is 1.  Now in this case we have
\[\mathbb{E}_{\mmu}\bbb{\exp\bb{\theta \sum_{j=1}^n f(s_j)}}\leq \exp\bb{\frac{\theta^2}{2}\cdot \alpha(\lambda)\cdot n},\]
as can be seen from our results with $d=1$.  We see that $n\alpha(\lambda)$ is the variance proxy for a sub-Gaussian random variable and so naturally upper bounds the asymptotic variance.  This lower bound is again attained for $P$ and $f$, and therefore the same holds for the matrix-valued function $F$. 

\subsection{Roadmap}
We begin with the quantity 
\[\mathbb{E}_{\mmu}\bbb{\norm{\prod_{j=1}^n \exp\bb{\frac{\theta e^{i\phi}}{2}F_j(s_j)}}_F^2}\]
for $\phi\in [-\pi/2, \pi/2]$ and $\theta>0$, developed first in \cite{sutter2017multivariate, garg2018}.  This will play the role of our moment generating function.  Using properties of the Kronecker product, in Section \ref{sec:starting} we show that this moment generating function is equal to
\[\inner{\mathbf{1}\otimes\operatorname{vec}(I_d)}{ E_1^{\theta/2}\bb{\prod_{j=1}^{n-1}E_j^{\theta/2}\tilde{P}E_{j+1}^{\theta/2}}E_n^{\theta/2}(\mathbf{1}\otimes\operatorname{vec}(I_d))}_{\mmu},\]
where $\tilde{P}$ is the lifted Markov operator and $E_j^{\theta/2}$ is a certain multiplication operator (Definition \ref{def:mul_operator}).  This quantity is a matrix version of a classical quantity that appears in the analysis of sums of Markov-dependent scalar random variables -- indeed, when $d=1$ these are equivalent.  Thus, we can apply a spectral analysis to bound this quantity.

In Section \ref{sec:operator} we proceed with this study; this section addresses the main challenges of extending to continuous state spaces.  Using Cauchy-Schwartz we bound the above by
\[\norm{E^{{\theta/2}^*}_1(\mathbf{1}\otimes\operatorname{vec}(I_d))}_{\mmu}\norm{E^{\theta/2}_n(\mathbf{1}\otimes\operatorname{vec}(I_d))}_{\mmu}\prod_{j=1}^{n-1}\norm{E_j^{\theta/2}\tilde{P}E_{j+1}^{\theta/2}}_{\mmu}.\]
Lemma \ref{lem:P_hat_props} is a symmetrizing result, allowing us to shift our attention to the lifted Leon-Perron operator $\hat{P}$ by giving the bound
\[\norm{E^{{\theta/2}^*}_1(\mathbf{1}\otimes\operatorname{vec}(I_d))}_{\mmu}\norm{E^{\theta/2}_n(\mathbf{1}\otimes\operatorname{vec}(I_d))}_{\mmu}\prod_{j=1}^{n-1}\norm{E_j^{\theta/2}\tilde{P}E_{j+1}^{\theta/2}}_{\mmu}\leq d\prod_{j=1}^n \norm{E_j^{\theta/2}\hat{P}E_j^{{\theta/2}^*}}_{\mmu}.\]
This lets us replace the operator $\tilde{P}$, which represents a chain that is not necessarily reversible, by $\hat{P}$, which represents a very simple reversible chain with many of the same important spectral properties as $\tilde{P}$.  In doing so, we obtain operators $E_j^{\theta/2}\hat{P}E_j^{{\theta/2}^*}$ that are now self-adjoint (and even positive semidefinite) on $\ell_2(\mmu\otimes \mathbf{1})$.  We focus on bounding the leading eigenvalue of each of these self-adjoint matrices.

To simplify the results, we then notice that the eigenvalues of $E_j^{\theta/2}\hat{P}E_j^{{\theta/2}^*}$ are equal to the eigenvalues of $E_{T_j}^{\theta/2}\hat{P}E_{T_j}^{\theta/2}$, for some other mulitiplication operator $E_{T_j}^{\theta/2}$ that is in fact \textit{real}, since the underlying matrix-valued functions $F_j$ are real-valued.  Our next main result states that the leading eigenvalue of $E_{T_j}^{\theta/2}\hat{P}E_{T_j}^{\theta/2}$ is a sort of limit for a corresponding moment generating function.  More precisely, we show that 
\[\lim_{n\to\infty}\frac{1}{n}\log \mathbb{E}_{\mmu}\bbb{\norm{\prod_{k=1}^n \exp\bb{\frac{\theta\cos(\phi)}{2}F_j(s_k)}}_F^2}=\log \norm{E_{T_j}^{\theta/2}\hat{P}E_{T_j}^{\theta/2}}_{\mmu},\]
where $s_1, \dots, s_n$ is driven by the Markov chain represented by the Leon-Perron operator $\hat{P}$; this result appears in Lemma \ref{lem:limit} and \ref{lem:simple_function}.  The realness assumption on the $F_j$ is important for the proof of this result -- we invoke an interesting generalization of the Perron-Frobenius theorem that applies to linear transformations leaving invariant a cone.  This result is our main technical innovation and is crucial in extending existing results to the continuous state space setting; conveniently it also allows us to sharpen known bounds in discrete state spaces.

Section \ref{sec:final} gives our final bounds.  We first show how to transfer bounds on our moment generating function to tail bounds; the technique is straightforward and uses the multi-matrix Golden-Thompson inequality from \cite{garg2018}.  We then give these tight bounds on the moment generating function under both Hoeffding and Bernstein-type assumptions.  For Hoeffding-type assumptions we use a coupling technique to exhibit a two-state chain that acts as the ``limit'' of our chain.  The leading eigenvalue of the corresponding operator can be solved exactly, giving optimal bounds.  For Bernstein-type assumptions, we use a robust linear algebraic approach to bound the operator norm of a related matrix directly; this approach will also give optimal results.

\section{Preliminaries}\label{sec:preliminaries}
\subsection{Notation}
Lower case, unbolded $a$ denotes scalars or scalar-valued functions, bold $\mathbf{a}$ denotes vectors or vector-valued functions, and upper case $A$ denotes matrices or matrix-valued functions.  

The operator $\otimes$ will denote the Kronecker product, where for matrices $A, B$ of size $a\times b$, $c\times d$, respectively, 
\[A\otimes B=\begin{bmatrix}
A_{1, 1}B & \dots & A_{1, b}B\\
\vdots & & \vdots \\
A_{a, 1}B & \dots & A_{a, b} B
\end{bmatrix}\]
of size $ac\times bd$.  This operation has an identification with the tensor product.  An important fact we use is $(A\otimes C)(B\otimes D)=AB\otimes CD$ for matrices of the appropriate sizes.  The Kronecker product has an interesting relationship with the vectorization operator $\operatorname{vec}$, where for a matrix $X$ of size $a\times b$, $\operatorname{vec}(X)$ is the flattened vector of size $ab$.  Importantly, $(B^\top \otimes A)\operatorname{vec}(X)=\operatorname{vec}(AXB)$, and so
\[\operatorname{tr}[AB]=\inner{\operatorname{vec}(I_d)}{\operatorname{vec}(AB)}=\inner{\operatorname{vec}(I_d)}{(B^\top \otimes A)\operatorname{vec}(I_d)}.\]

\subsection{Markov chains}
Let $\mathcal{X}$ be a general state space with  $\sigma$-algebra $\mathcal{B}:=\mathcal{B}(\mathcal{X})$.  Let $P$ be a Markov kernel on $\mathcal{X}$ defined, for a sequence of random variables $X_1, \dots, X_n$, in the standard way as
\[P(x, B)=\operatorname{Pr}\bb{X_k\in B\mid X_{k-1}=x},\quad \forall B\in\mathcal{B}\]
with stationary measure $\mmu$ so that
\[\mmu(B)=\int P(x, B)\mmu(dx), \quad \forall B\in\mathcal{B}.\]
Define $\ell_2(\mmu)$ as 
\[\ell_2(\mmu)=\{h:\mathcal{X}\to\bbC\mid \mathbb{E}_{\mmu}[h(x)^2]<\infty\}\]
for $h:\mathcal{X}\to\bbC$ measurable.  This is a a Hilbert space equipped with the following inner product:
\[\inner{f}{g}_{\mmu}=\int f(x)g(x)\, d\mmu(x), \quad \forall f, g\in \ell_2(\mmu)\]
and corresponding norm $\norm{f}_{\mmu}=\sqrt{\inner{f}{f}_{\mmu}}$.  The norm of a linear operator $T$ is defined in the usual way:
\[\norm{T}_{\mmu}=\sup_{\norm{h}_{\mmu}=1}\norm{Th}_{\mmu}.\]

A transition kernel $P$ acts as an operator on $\ell_2(\mmu)$:
\[(Ph)(x)=\int h(y)P(x, dy), \quad \forall x\in \mathcal{X}, h\in \ell_2(\mmu).\]

The projection operator $\Pi$ corresponding to the distribution $\mmu$ is defined as $(\Pi h)(x)=\mathbb{E}_{\mmu}[h]\mathbf{1}$ ($\mathbf{1}$ is the function such that $\mathbf{1}(x)=1$ for all $x$); this is a rank-1 operator by definition, and is indeed a projection onto $\mathbf{1}$ as $\mathbb{E}_{\mmu}[h]=\inner{\mathbf{1}}{h}_{\mmu}$.  If $\mmu$ is stationary for the transition kernel $P$ then $P\Pi=\Pi P=\Pi$.

An important subset of elements of $\ell_{2}(\mmu)$ will be the class of ``mean-zero functions,'' denoted
\[\ell_2^0(\mmu)=\{h\in\ell_2(\mmu)\mid \Pi h=0\}.\]
Then we have the following definition:
\begin{definition}[Absolute spectral gap]\label{def:absolute_spectral}
    A Markov kernel $P$ with stationary measure $\mmu$ admits an absolute spectral gap $1-\lambda(P)$ if 
    \[\lambda(P):=\sup_{h\in\ell_2^0(\mmu), h\neq 0}\frac{\norm{Ph}_{\mmu}}{\norm{h}_{\mmu}}=\norm{P-\Pi}_{\mmu}<1.\]
\end{definition}
Note that $\lambda(P)\leq 1$ always, as $\norm{Ph}_{\mmu}\leq \norm{h}_{\mmu}$ by Jensen's inequality, with equality for $h=\mathbf{1}$.  When it is clear from context we will just use $\lambda=\lambda(P)$.  The absolute spectral gap characterizes the convergence of a Markov chain to its invariant measure \cite{Rudolf2012}.  For reversible finite-state chains, the value $\lambda(P)$ corresponds to the second largest eigenvalue and the existence of the gap corresponds to ergodicity.  A Markov chain driven by $P$ is reversible if and only if $P$ is a self-adjoint operator on $\ell_2(\mmu)$.  

We define what is known as a Leon-Perron operator \cite{leon2004, fan2021, jiang2018} that will be important in the sequel:
\begin{definition}[Leon-Perron operator]
    Let $P$ be a Markov kernel with stationary measure $\mmu$.  For a constant $c\in [0, 1]$, define $\hat{P}_c$ as
    \[\hat{P}_c:= cI+(1-c)\Pi.\]
    If $c=\lambda(P)$ we drop the subscript and call $\hat{P}$ the Leon-Perron version of $P$.
\end{definition}
We can interpret $\hat{P}_c$ as a transition kernel such that at a state, it stays at that state with probability $c$ or samples a new state independently from $\mmu$ with probability $1-c$.  Note that if $P$ admits an absolute spectral gap then $\hat{P}$ does so as well with the same absolute spectral gap.  

\subsection{Behavior in tensored space}
As we work with matrices, the above will have to be lifted to a tensored space as necessary.  We first consider the product space $\mathbb{C}^{d}\otimes \mathbb{C}^d\simeq \mathbb{C}^{d^2}$ with inner product induced from the standard Euclidean product, i.e., $\inner{\mathbf{a}\otimes \mathbf{c}}{\mathbf{b}\otimes \mathbf{d}}=\inner{\mathbf{a}}{\mathbf{b}}\inner{\mathbf{c}}{\mathbf{d}}$, and extend by linearity.  Oftentimes we will simply identify an element as $\mathbf{y}\in\mathbb{C}^{d^2}$; it is straightforward that the standard Euclidean inner product on this larger space is equivalent to the tensored inner product.  So it is no loss to move between one representation to the other.

Now define $\ell_2(\mmu\otimes \mathbf{1})$ as the following ``lift'' of $\ell_2(\mmu)$: formally, define
\[\ell_2(\mmu\otimes \mathbf{1})=\{\mathbf{h}:\mathcal{X}\to\mathbb{C}^{d^2}\mid \mathbb{E}_{\mmu}[\norm{\mathbf{h}(x)}^2]<\infty\}\]
for $\mathbf{h}:\mathcal{X}\to\bbC^{d^2}$ measurable as a vector-valued function.  This is equipped with the inner product
\[\inner{\mathbf{f}}{\mathbf{g}}_{\mmu}=\int \inner{\mathbf{f}(x)}{\mathbf{g}(x)}\,d\mmu(x)\]
and corresponding norm $\norm{\mathbf{f}}_{\mmu}=\sqrt{\inner{\mathbf{f}}{\mathbf{f}}_{\mmu}}$.  This gives the understanding of $\ell_2(\mmu\otimes \mathbf{1})$ as a direct integral of Hilbert spaces indexed by the points of $\mathcal{X}$; decomposable elements in this space are understood to be of the form $f\otimes \mathbf{v}$ for $f\in\ell_2(\mmu), \mathbf{v}\in\mathbb{C}^{d^2}$ such that $(f\otimes \mathbf{v})(x)=f(x)\mathbf{v}$.   Norms of linear operators are thus defined with respect to this norm in the usual way.

For a Markov operator $P$ on $\ell_2(\mmu)$ we can define the operator $\tilde{P}:=P\otimes I_{d^2}$ as an operator on $\ell_2(\mmu\otimes \mathbf{1})$.  This operator acts as 
\[(\tilde{P}\mathbf{h})(x)=\int \mathbf{h}(x)P(x, dy), \quad \forall x\in\mathcal{X}, \mathbf{h}\in\ell_2(\mmu\otimes \mathbf{1});\]
note that $(\tilde{P}\mathbf{h})(x)$ is a vector.  Whereas for the operator $P$ the leading eigenfunction was $\mathbf{1}$ and spans the entire eigenspace for the leading eigenvalue of $1$ (assuming nonzero spectral gap), for $\tilde{P}$ the leading eigenspace is $d^2$-dimensional, spanned by $\{\mathbf{1}\otimes \mathbf{e}_1, \dots, \mathbf{1}\otimes \mathbf{e}_{d^2}\}$.  This eigenspace is $\mathbf{1}\otimes \mathbb{C}^{d^2}$.

For a decomposable element $f\otimes \mathbf{v}$, where $f\in \ell_2(\mmu), \mathbf{v}\in\mathbb{C}^{d^2}$, projection onto $\mathbf{1}\otimes \mathbb{C}^{d^2}$ is $\inner{\mathbf{1}}{f}_{\mmu}(\mathbf{1}\otimes \mathbf{v})=\mathbb{E}_{\mmu}[f](\mathbf{1}\otimes \mathbf{v})$.  For the projection operator $\Pi$ we can define the lifted version as $\tilde{\Pi}:=\Pi\otimes I_{d^2}$, again satisfying $\tilde{\Pi}\tilde{P}=\tilde{P}\tilde{\Pi}=\tilde{\Pi}$, that performs exactly this operation, i.e., $\tilde{\Pi}(f\otimes \mathbf{v})=\inner{\mathbf{1}}{f}_{\mmu}(\mathbf{1}\otimes \mathbf{v})$.  This can all be extended by linearity and convergence in $\ell_2$.  Then there is a corresponding lift of ``mean-zero functions'' as
\[\ell_2^0(\mmu\otimes \mathbf{1})=\{\mathbf{h}\in \ell_2(\mmu\otimes \mathbf{1})\mid \tilde{\Pi}\mathbf{h}=0\},\]
and can also be identified as $\ell_2^0(\mmu\otimes \mathbf{1})=\ell_2^0(\mmu)\otimes \mathbb{C}^{d^2}$.  Similarly, define
\[\lambda(\tilde{P})=\sup_{\mathbf{h}\in \ell_2^0(\mmu\otimes \mathbf{1}), \mathbf{h}\neq 0}\frac{\norm{\tilde{P}\mathbf{h}}_{\mmu}}{\norm{\mathbf{h}}_{\mmu}}=\norm{\tilde{P}-\tilde{\Pi}}_{\mmu}.\]
The following lemma, Lemma 3 from \cite{qiu2020}, will imply that if $P$ admits an absolute spectral gap, then so does $\tilde{P}$ with the same absolute spectral gap.  Thus, many of the essential spectral properties of $P$ and $\tilde{P}$ are the same.  The proof is deferred to the appendix.
\begin{lemma}\label{lem:lambda_Phat}
    $\lambda(P)=\lambda(\tilde{P})$.
\end{lemma}

Lastly, we can define Leon-Perron operators of $\tilde{P}$ as the corresponding lifts of Leon-Perron operators of $P$, which is tensorization by $I_{d^2}$.  With some abuse of notation, we will often interchange the notation $\hat{P}$ to refer to both Leon-Perron operators of $P$ and Leon-Perron operators of $\tilde{P}$, but the usage will be clear from context.

Recall for a scalar-valued function $g$, the multiplication operator $M_g$ is defined as $(M_g h)(x)=g(x)h(x)$.  For a matrix-valued function $G$ we have the following definition:
\begin{definition}[Multiplication operator]\label{def:mul_operator}
    Let $G$ be a matrix-valued function mapping $\mathcal{X}$ to $\bbC^{d^2\times d^2}$ matrices.  Define $M_G$ be the operator on $\ell_2(\mmu\otimes \mathbf{1})$ such that for any $\mathbf{h}\in \ell_2(\mmu\otimes \mathbf{1})$, $(M_G \mathbf{h})(x)=G(x)\mathbf{h}(x)$.  We also call these \textit{block diagonal} operators.
\end{definition}
In other words, $M_G$ is multiplication by a function $G$ in the direct integral of Hilbert spaces \cite{nauimark1964}.  We will often use the multiplication operator $E^{\theta}_H$ defined as $(E^{\theta}\mathbf{h})(x)=\exp(\theta H(x))\mathbf{h}(x)$, sometimes dropping the subscript $H$ when it is clear from context.

\section{Starting point}\label{sec:starting}
Our starting point is the following matrix moment generating function:
\begin{equation}\label{eq:matrix_mgf}
    \mathbb{E}_{\mmu}\bbb{\norm{\prod_{j=1}^n \exp\bb{\frac{\theta e^{i\phi}}{2}F_j(s_j)}}_F^2},
\end{equation}
where $\phi\in [-\pi/2, \pi/2]$, $\theta>0$, and $s_1, \dots, s_n$ is a sequence of states driven by the Markov chain with transition matrix $P$.

To simplify this moment generating function, we first use the property that $\operatorname{tr}[AB^\top]=\operatorname{vec} (I_d)^\top (A\otimes B)\operatorname{vec}(I_d)$, where $\otimes$ is the Kronecker product and $I_d$ is the $d\times d$ identity matrix to write
\begin{equation}
\begin{aligned}
& \norm{\prod_{j=1}^n \exp\bb{\frac{\theta e^{i\phi}}{2}F_j(s_j)}}_F^2\\
=\, & \operatorname{tr}\bbb{\prod_{j=1}^n \exp\bb{\frac{\theta e^{i\phi}}{2}F_j(s_j)}\prod_{j=n}^1\exp\bb{\frac{\theta e^{-i\phi}}{2}F_j(s_j)}}\\
=\, & \operatorname{vec}(I_d)^\top\bb{\prod_{j=1}^n \exp\bb{\frac{\theta e^{i\phi}}{2}F_j(s_j)}\otimes \prod_{j=1}^n \exp\bb{\frac{\theta e^{-i\phi}}{2}F_j(s_j)}}\operatorname{vec}(I_d)
\end{aligned}
\end{equation}
and then use successive applications of the property $(AC)\otimes (BD)=(A\otimes B)(C\otimes D)$ to write the above as
\[
    \operatorname{vec}(I_d)^\top\bb{\prod_{j=1}^n \exp\bb{\frac{\theta e^{i\phi}}{2}F_j(s_j)}\otimes \exp\bb{\frac{\theta e^{-i\phi}}{2}F_j(s_j)}}\operatorname{vec}(I_d).
\]
Lastly, we use the fact that $\exp(A)\otimes \exp(B)=\exp(A\otimes I_d+I_d\otimes B)$ to ultimately write
\begin{equation}
    \operatorname{tr}\bbb{\prod_{j=1}^n \exp\bb{\frac{\theta e^{i\phi}}{2}F_j(s_j)}\prod_{j=n}^1\exp\bb{\frac{\theta e^{-i\phi}}{2}F_j(s_j)}}=\operatorname{vec}(I_d)^* \bb{\prod_{j=1}^n \exp(\theta H_j(s_j))}\operatorname{vec}(I_d)
\end{equation}
where $H_1, \dots, H_n$ is a sequence of matrix-valued function on $\mathcal{X}$, implicitly with respect to $\phi$, defined as
\begin{equation}
    H_j(x)=\frac{e^{i\phi}}{2}F_j(x)\otimes I_d+\frac{ e^{-i\phi}}{2}I_d\otimes F_j(x)\in \mathbb{C}^{d^2\times d^2}
\end{equation}
for $x\in\mathcal{X}$.  Note that unfortunately $H_j(x)$ is not in general Hermitian for any $x$, though its real and imaginary parts are both symmetric.  We emphasize the identity
\[\exp(A)\otimes \exp(B)=\exp(A\otimes I_d+I_d\otimes B)\]
as we will switch back and forth from these two forms as needed.  

The above allows us to focus on the matrices $H_j$, which are measurable and satisfy many of the same probabilistic properties as $F_j$ -- a more precise statement is given in the next section.  We then write
\begin{equation}\label{eq:seven}
    \begin{aligned}
    & \mathbb{E}_{\mmu}\bbb{\prod_{j=1}^n \exp(\theta H_j(s_j))}\\
    =\, & \int P(s_1, ds_2)\dots P(s_{n-1}, ds_n)\prod_{j=1}^n \exp(\theta H_j(s_j))\,d\mmu(s_1)\\
    =\, & \int d\mmu(s_1)\int P(s_1, ds_2)\exp(\theta H_1(s_1))\dots \int P(s_{n-1}, ds_n) \exp(\theta H_{n-1}(s_{n-1}))\exp(\theta H_n(s_n)).
    \end{aligned}
\end{equation}

Recall that $E^{\theta}_j$ the multiplication operator defined as $(E^{\theta}_j\mathbf{h})(x)=\exp(\theta H_j(x))\mathbf{h}(x)$ for any vector-valued function $\mathbf{h}$ -- see Definition \ref{def:mul_operator}, and recall that $\tilde{P}=P\otimes I_{d^2}$.  $E_j^\theta$ and $\tilde{P}$ act via
\[(\tilde{P}E_{j}^\theta\mathbf{h})(x)=\int \exp(\theta H_j(y))\mathbf{h}(y)P(x, dy)\]
and 
\[(E_{j}^{\theta/2}\tilde{P}E_{j}^{\theta/2}\mathbf{h})(x)=\int \exp\bb{\frac{\theta}{2} H_j(x)}\exp\bb{\frac{\theta}{2} H_j(y)}\mathbf{h}(y)P(x, dy)\]
as operators.  Then simplifying Equation \ref{eq:seven}, and using the definition of the inner product on $\ell_2(\mmu\otimes \mathbf{1})$, we have
\begin{equation}\label{eq:markov_matrix_mgf}
    \mathbb{E}_{\mmu}\bbb{\prod_{j=1}^n \exp(\theta H_j(s_j))}=\inner{\mathbf{1}\otimes \operatorname{vec}(I_d)}{E^{\theta/2}_1\bb{\prod_{j=1}^{n-1} E^{\theta/2}_j\tilde{P}E^{\theta/2}_{j+1}}E^{\theta/2}_n(\mathbf{1}\otimes \operatorname{vec}(I_d))}_{\mmu}.
\end{equation}

\section{Bounding the operator norm}\label{sec:operator}
Having derived the above, in this section we give a series of bounds that shift our focus to bounding the norm of a perturbed operator; we address the main technical challenges of studying general continuous state spaces. 

From Equation \ref{eq:markov_matrix_mgf}, we apply Cauchy-Schwartz and submultiplicativity of the spectral norm to obtain
\begin{equation}\label{eq:mgf_step1}
    \begin{aligned}
    & \inner{\mathbf{1}\otimes \operatorname{vec}(I_d)}{E^{\theta/2}_1\bb{\prod_{j=1}^{n-1} E^{\theta/2}_j\tilde{P}E^{\theta/2}_{j+1}}E^{\theta/2}_n(\mathbf{1}\otimes \operatorname{vec}(I_d))}_{\mmu}\\
    =\, & \inner{E^{{\theta/2}^*}_1 (\mathbf{1}\otimes \operatorname{vec}(I_d))}{\bb{\prod_{j=1}^{n-1} E^{\theta/2}_j\tilde{P}E^{\theta/2}_{j+1}}E^{\theta/2}_n(\mathbf{1}\otimes \operatorname{vec}(I_d))}_{\mmu}\\
    \leq\, &  \norm{E^{{\theta/2}^*}_1(\mathbf{1}\otimes\operatorname{vec}(I_d))}_{\mmu}\norm{E^{\theta/2}_n(\mathbf{1}\otimes\operatorname{vec}(I_d))}_{\mmu}\prod_{j=1}^{n-1}\norm{E_j^{\theta/2}\tilde{P}E_{j+1}^{\theta/2}}_{\mmu}.
    \end{aligned}
\end{equation}
Now $\tilde{P}$ possesses many of the same spectral properties as $P$, most notably the result from Lemma \ref{lem:lambda_Phat}.  Recall that $\hat{P}$ is the ``lifted'' Leon-Perron operator for $\tilde{P}$, defined as $\hat{P}=(\lambda I+(1-\lambda)\Pi)\otimes I_{d^2}=\lambda I+(1-\lambda)\tilde{\Pi}$ with $\tilde{\Pi}=\Pi\otimes I_{d^2}$.  The following lemma allows us to replace $\tilde{P}$ with its Leon-Perron version $\hat{P}$.
\begin{lemma}\label{lem:P_hat_props}
The operator $\tilde{P}$ and its Leon-Perron version $\hat{P}$ satisfy the following:
\begin{enumerate}
    \item For any $\mathbf{g}, \mathbf{h}\in \ell_2(\mmu\otimes \mathbf{1})$, $\abs{\inner{\mathbf{g}}{\tilde{P} \mathbf{h}}_{\mmu}}\leq \inner{\mathbf{g}}{\hat{P}\mathbf{g}}_{\mmu}^{1/2}\inner{ \mathbf{h}}{\hat{P} \mathbf{h}}_{\mmu}^{1/2}$.
    \item For operators $S_1, S_2$ on $\ell_2(\mmu\otimes \mathbf{1})$, $\norm{S_1\tilde{P}S_2}_{\mmu}\leq \norm{S_1\hat{P}S_1^*}_{\mmu}^{1/2}\norm{S_2\hat{P}S_2^*}_{\mmu}^{1/2}$.
    \item For any multiplication operator $M_G$ with respect to a measurable matrix valued function $G$, and for any vector $\mathbf{v}$, $\norm{M_G(\mathbf{1}\otimes \mathbf{v})}_{\mmu}\leq \norm{\mathbf{v}}\norm{M_G\hat{P}M_G^*}_{\mmu}^{1/2}$.
\end{enumerate}
\end{lemma}
\begin{proof} 
For (1), we have
    \begin{align*}
        \abs{\inner{\mathbf{g}}{\tilde{P} \mathbf{h}}_{\mmu}}
        &=\abs{\inner{\mathbf{g}}{(\tilde{P}-\tilde{\Pi}+\tilde{\Pi})\mathbf{h}}}\\
        &=\abs{\inner{\mathbf{g}}{\tilde{\Pi}\mathbf{h}}+\inner{\mathbf{g}}{(\tilde{P}-\tilde{\Pi})\mathbf{h}}}\\
        &\leq \abs{\inner{\mathbf{g}}{\tilde{\Pi} \mathbf{h}}_{\mmu}}+ \abs{\inner{\mathbf{g}-\tilde{\Pi}\mathbf{g}}{(\tilde{P}-\tilde{\Pi})( \mathbf{h}-\tilde{\Pi} \mathbf{h})}_{\mmu}}\\
        &\leq \abs{\inner{\mathbf{g}}{\tilde{\Pi} \mathbf{h}}_{\mmu}}+ \lambda(\tilde{P})\norm{(I-\tilde{\Pi})\mathbf{g}}_{\mmu}\norm{(I-\tilde{\Pi})\mathbf{h}}_{\mmu}.
    \end{align*}
    where the first inequality follows because $\tilde{\Pi}\tilde{P}=\tilde{P}\tilde{\Pi}=\tilde{\Pi}$ and because $\tilde{\Pi}$ is a projection and is thus self-adjoint on $\ell_2(\mmu\otimes \mathbf{1})$.  Because of this we also have $\inner{\mathbf{g}}{\tilde{\Pi}\mathbf{h}}_{\mmu}=\inner{\tilde{\Pi}\mathbf{g}}{\tilde{\Pi}\mathbf{h}}_{\mmu}$.  Therefore, by Cauchy-Schwartz, $\abs{\inner{\mathbf{g}}{\tilde{\Pi}\mathbf{h}}_{\mmu}}\leq \norm{\tilde{\Pi}\mathbf{g}}_{\mmu}\norm{\tilde{\Pi}\mathbf{h}}_{\mmu}$.  Then
    \begin{align*}
        \abs{\inner{\mathbf{g}}{\tilde{P}\mathbf{h}}_{\mmu}}
        &\leq \abs{\inner{\mathbf{g}}{\tilde{\Pi}\mathbf{h}}_{\mmu}}+ \lambda(\tilde{P})\norm{(I-\tilde{\Pi})\mathbf{g}}_{\mmu}\norm{(I-\tilde{\Pi})\mathbf{h}}_{\mmu}\\
        &\leq \norm{\tilde{\Pi}\mathbf{g}}_{\mmu}\norm{\tilde{\Pi}\mathbf{h}}_{\mmu}+\lambda(\tilde{P})\norm{(I-\tilde{\Pi})\mathbf{g}}_{\mmu}\norm{(I-\tilde{\Pi})\mathbf{h}}_{\mmu}\\
        &\leq \sqrt{\lambda \norm{(I-\tilde{\Pi})\mathbf{g}}_{\mmu}^2+\norm{\tilde{\Pi}\mathbf{g}}_{\mmu}^2}\cdot \sqrt{\lambda \norm{(I-\tilde{\Pi})\mathbf{h}}_{\mmu}^2+\norm{\tilde{\Pi}\mathbf{h}}_{\mmu}^2}\\
        &=\inner{\mathbf{g}}{\hat{P}\mathbf{g}}_{\mmu}^{1/2}\inner{\mathbf{h}}{\hat{P}\mathbf{h}}_{\mmu}^{1/2}.
    \end{align*}
    
For (2), we have
    \begin{align*}
        \norm{S_1\tilde{P} S_2}_{\mmu}
        &=\sup_{\norm{\mathbf{g}}_{\mmu}=\norm{\mathbf{h}}_{\mmu}=1}\abs{\inner{\mathbf{g}}{S_1\tilde{P}S_2 \mathbf{h}}_{\mmu}}\\
        &=\sup_{\norm{\mathbf{g}}_{\mmu}=\norm{\mathbf{h}}_{\mmu}=1}\abs{\inner{S_1^* \mathbf{g}}{\tilde{P}S_2\mathbf{h}}_{\mmu}}\\
        &\leq \sup_{\norm{\mathbf{g}}_{\mmu}=\norm{\mathbf{h}}_{\mmu}=1} \inner{S_1^* \mathbf{g}}{\hat{P}S_1^* \mathbf{g}}_{\mmu}^{1/2}\inner{S_2\mathbf{h}}{\hat{P}S_2\mathbf{h}}_{\mmu}^{1/2}\\
        &=\norm{S_1\hat{P}S_1^*}_{\mmu}^{1/2}\norm{S_2^* \hat{P}S_2}_{\mmu}^{1/2}.
    \end{align*}
    The result follows by noting that $\norm{S_1^* \hat{P}S_1}_{\mmu}=\norm{S_1\hat{P}S_1^*}_{\mmu}$.
    
For (3), the case $G=0$ is clear.  Otherwise, 
    \begin{align*}
        \norm{M_G\hat{P}M_G^*}_{\mmu}
        &\geq \frac{\inner{M_G(\mathbf{1}\otimes \mathbf{v})}{M_G\hat{P}M_G^* M_G (\mathbf{1}\otimes \mathbf{v})}_{\mmu}}{\norm{M_G(\mathbf{1}\otimes \mathbf{v})}^2_{\mmu}}\\
        &= \frac{\inner{M_G^*M_G(\mathbf{1}\otimes \mathbf{v})}{\tilde{\Pi}M_G^* M_G (\mathbf{1}\otimes \mathbf{v})}_{\mmu}}{\norm{M_G(\mathbf{1}\otimes \mathbf{v})}^2_{\mmu}}+ \frac{\lambda\inner{M_G^*M_G(\mathbf{1}\otimes \mathbf{v})}{(I-\tilde{\Pi})M_G^* M_G (\mathbf{1}\otimes \mathbf{v})}_{\mmu}}{\norm{M_G(\mathbf{1}\otimes \mathbf{v})}^2_{\mmu}}\\
        &=\frac{\inner{M_G^*M_G(\mathbf{1}\otimes \mathbf{v})}{\tilde{\Pi}M_G^* M_G (\mathbf{1}\otimes \mathbf{v})}_{\mmu}}{\norm{M_G(\mathbf{1}\otimes \mathbf{v})}^2_{\mmu}}+ \frac{\lambda\norm{(I-\tilde{\Pi})M_G^* M_G(\mathbf{1}\otimes \mathbf{v})}^*_{\mmu}}{\norm{M_G(\mathbf{1}\otimes \mathbf{v})}^2_{\mmu}}\\
        &\geq \frac{\inner{M_G^*M_G(\mathbf{1}\otimes \mathbf{v})}{\tilde{\Pi}M_G^* M_G (\mathbf{1}\otimes \mathbf{v})}_{\mmu}}{\norm{M_G(\mathbf{1}\otimes \mathbf{v})}^2_{\mmu}},
    \end{align*}
    where the second to last line follows because $I-\tilde{\Pi}$ is a projection.  We can express these terms more explicitly.  First, for the denominator, we have $(M_G(\mathbf{1}\otimes \mathbf{v}))(x)=G(x)\mathbf{v}$.  Therefore, we have 
    \[\norm{M_G(\mathbf{1}\otimes \mathbf{v})}_{\mmu}^2=\inner{\mathbf{v}}{\mathbb{E}_{\mmu}[G(x)^*G(x)]\mathbf{v}}.\]
    Similarly, we have $\tilde{\Pi}M_G^* M_G(\mathbf{1}\otimes \mathbf{v})(x)=\mathbb{E}_{\mmu}[G(y)^*G(y)]\mathbf{v}$ for any $x\in\mathcal{X}$.  Then
    \[\inner{M_G^* M_G(\mathbf{1}\otimes \mathbf{v})}{\tilde{\Pi}M_G^* M_G (\mathbf{1}\otimes \mathbf{v})}_{\mmu}=\inner{\mathbb{E}_{\mmu}[G(x)^*G(x)]\mathbf{v}}{\mathbb{E}_{\mmu}[G(x)^*G(x)]\mathbf{v}}.\]
    Now for any Hermitian matrix $A$ and any vector $\mathbf{v}\neq 0$ it holds that $\inner{A\mathbf{v}}{A\mathbf{v}}\geq \inner{\mathbf{v}}{A\mathbf{v}}^2/\norm{\mathbf{v}}^2$ as a simple consequence of Cauchy-Schwartz.  Therefore, if $A$ is positive semidefinite, then $\inner{A\mathbf{v}}{A\mathbf{v}}/\inner{\mathbf{v}}{A\mathbf{v}}\geq \inner{\mathbf{v}}{A\mathbf{v}}/\norm{\mathbf{v}}^2$.  As $\mathbb{E}_{\mmu}[G(x)^*G(x)]$ is indeed positive semidefinite,
    \begin{align*}
        \norm{M_G\hat{P}M_G^*}_{\mmu}
        &\geq \frac{\inner{M_G^* M_G(\mathbf{1}\otimes \mathbf{v})}{\tilde{\Pi}M_G^* M_G (\mathbf{1}\otimes \mathbf{v})}_{\mmu}}{\norm{M_G(\mathbf{1}\otimes \mathbf{v})}^2_{\mmu}}\\
        &=\frac{\inner{\mathbb{E}_{\mmu}[G(x)^*G(x)]\mathbf{v}}{\mathbb{E}_{\mmu}[G(x)^*G(x)]\mathbf{v}}}{\inner{\mathbf{v}}{\mathbb{E}_{\mmu}[G(x)^*G(x)]\mathbf{v}}}\\
        &\geq \frac{\inner{\mathbf{v}}{\mathbb{E}_{\mmu}[G(x)^*G(x)]\mathbf{v}}}{\norm{\mathbf{v}}^2}\\
        &=\frac{\norm{M_G(\mathbf{1}\otimes \mathbf{v})}^2_{\mmu}}{\norm{\mathbf{v}}^2}.
    \end{align*}
    Rearranging finishes the proof.
\end{proof}
Going back to Equation \ref{eq:mgf_step1}, applying Lemma \ref{lem:P_hat_props} parts (2) and (3), and noticing that $\norm{\operatorname{vec}(I_d)}=\sqrt{d}$, we establish
\begin{equation}\label{eq:operator_bound}
    \begin{aligned}
        & \inner{\mathbf{1}\otimes \operatorname{vec}(I_d)}{E^{\theta/2}_1\bb{\prod_{j=1}^{n-1} E^{\theta/2}_j\tilde{P}E^{\theta/2}_{j+1}}E^{\theta/2}_n(\mathbf{1}\otimes \operatorname{vec}(I_d))}_{\mmu}\\
        \leq \, & \norm{E^{{\theta/2}^*}_1(\mathbf{1}\otimes\operatorname{vec}(I_d))}_{\mmu}\norm{E^{\theta/2}_n(\mathbf{1}\otimes\operatorname{vec}(I_d))}_{\mmu}\prod_{j=1}^{n-1}\norm{E_j^{\theta/2}\tilde{P}E_{j+1}^{\theta/2}}_{\mmu}\\
        \leq \, & d \prod_{j=1}^n \norm{E_j^{\theta/2}\hat{P}E_j^{{\theta/2}^*}}_{\mmu}.
    \end{aligned}
\end{equation}
Thus, we have transferred the study of the problem from $\tilde{P}$, which represents a general, nonreversible Markov chain, to $\hat{P}$, which represents a reversible chain and thus is much simpler to analyze.  Therefore, the operators $E_j^{\theta/2}\hat{P}E_j^{{\theta/2}^*}$ are self-adjoint on $\ell_2(\mmu\otimes \mathbf{1})$.  Now we focus on bounding the leading eigenvalues of these operators; we make one more simplification:
\begin{proposition}\label{prop:T_operator}
    Let $T_j(x)=\frac{\cos(\phi)}{2}(F_j(x)\otimes I_d+I_d\otimes F(x)))$.  Let $E_{T_j}^\theta$ be the operator defined as $(E_{T_j}^\theta \mathbf{h})(x)=\exp(\theta T_j(x))\mathbf{h}(x)$.  Then the norms of $E_j^{\theta/2}\hat{P}E_j^{{\theta/2}^*}$ and $E_{T_j}^{\theta/2}\hat{P}E_{T_j}^{\theta/2}$ are the same.
\end{proposition}
\begin{proof}
    The operator $E^{\theta/2}_j\hat{P}E_j^{{\theta/2}^*}$ is similar to the operator $\hat{P}E_j^{{\theta/2}^*}E^{\theta/2}_j$.  The operator $E_j^{{\theta/2}^*} E_j^{\theta/2}$ acts via
    \begin{align*}
        (E_j^{{\theta/2}^*}E_j^{\theta/2}\mathbf{h})(x)
        &=\exp\bb{\frac{\theta}{2}(H_j^*(x)+H_j(x))}\mathbf{h}(x)\\
        &=\exp\bb{\frac{\theta}{2}\bb{\frac{e^{-i\phi}+e^{i\phi}}{2}F_j(x)\otimes I_d+\frac{e^{i\phi}+e^{-i\phi}}{2}I_d\otimes F_j(x)}}\mathbf{h}(x)\\
        &=\exp\bb{\frac{\theta\cos(\phi)}{2}(F_j(x)\otimes I_d+I_d\otimes F_j(x))}\mathbf{h}(x)\\
        &=(E_{T_j}^{\theta}\mathbf{h})(x),
    \end{align*}
    where the first equality holds because $H_j(x), H_j^*(x)$ commute.  As $\hat{P}E_{T_j}^{\theta}$ is similar to $E_{T_j}^{\theta/2}\hat{P}E_{T_j}^{\theta/2}$, we see that $E_j^{\theta/2}\hat{P}E_j^{{\theta/2}^*}$ and $E_{T_j}^{\theta/2}\hat{P}E_{T_j}^{\theta/2}$ have the same spectrum; as they are both self-adjoint, the spectral radius equals the operator norm.
\end{proof}
At this point, we have reduced our problem to a considerably simpler form.  We can now show that the operators $T_j$ satisfy many of the same probabilistic properties as $F_j$; the proof is in the appendix.
\begin{proposition}\label{prop:T_properties}
    Let $T$ be the operator defined as $T(x)=\frac{\cos(\phi)}{2}(F(x)\otimes I_d+I_d\otimes F(x))$ with $\phi\in [-\pi/2, \pi/2]$ so that $\cos(\phi)$ is always nonnegative.  Then
    \begin{enumerate}
        \item If $\mathbb{E}[F(x)]=0$, then $\mathbb{E}[T(x)]=0$,
        \item If $aI\preceq F(x)\preceq bI$ for all $x\in\mathcal{X}$, then $a\cos(\phi)I\preceq T(x)\preceq b\cos(\phi)I$ for all $x\in\mathcal{X}$,
        \item If $\norm{\mathbb{E}[F(x)^2]}\leq \mathcal{V}$, then $\norm{\mathbb{E}[T(x)^2]}\leq \cos^2(\phi)\mathcal{V}$,
    \end{enumerate}
    where all expectations are taken with respect to a measure on $\mathcal{X}$, and $F$ and $T$ are measurable.
\end{proposition}

\subsection{The leading eigenvalue as the limit}
In this subsection we will prove a limit lemma, justifying the use of the leading eigenvalue as a characterization of the moment generating function, mirroring the scalar case as laid out in \cite{dembo2009} -- essentially, what we will prove is that the log limit of the moment generating function is the logarithm of the leading eigenvalue.  This will allow us to shift between time-dependent and time-independent functions as needed.  The main difficulty in extending this theory from the scalar case is proving the limit statement in part (3) of Lemma \ref{lem:simple_function}, and in particular, a specific lower bound in the proof.  This lower bound is not difficult to prove in the scalar setting; however, the matrix setting is more subtle and requires some more care.

Our main lemma is the following:
\begin{lemma}\label{lem:limit}
    Let $F$ be a map from $\mathcal{X}$ to real symmetric $d\times d$ matrices with $aI\preceq F(x)\preceq bI$ for all $x\in\mathcal{X}$.  Let $T(x)=\frac{1}{2}(F(x)\otimes I_d+I_d\otimes F(x))$, $x\in \mathcal{X}$, and let $E_T^{\theta/2}$ be the operator defined as $(E_T^{\theta}\mathbf{h})(x)=\exp(\theta T(x))\mathbf{h}(x)$.  Let $\mmu$ be a distribution, let $s_1, \dots, s_n$ be a sequence of states driven by the stationary Markov chain with transition matrix $\lambda I+(1-\lambda)\mathbf{1}\mmu^\top$, and let $\hat{P}$ be this transition matrix tensored with $I_{d^2}$.  Then
    \[\lim_{n\to\infty}\frac{1}{n}\log \mathbb{E}_{\mmu}\bbb{\norm{\prod_{j=1}^n \exp\bb{\frac{\theta }{2}F(s_j)}}_F^2}= \log\norm{E^{\theta/2}_T\hat{P}E^{\theta/2}_T}_{\mmu}. \]
\end{lemma}

In other words, we take $s_1, \dots, s_n$ to be driven by a Leon-Perron operator, where at a given state the chain stays at that state with probability $\lambda$ and samples a state from $\mmu$ with probability $1-\lambda$.  We first address the case of a simple matrix-valued function $F$.

\begin{lemma}\label{lem:simple_function}
    Let $\hat{P}=(\lambda I+(1-\lambda)\Pi)\otimes I_{d^2}$ be a Leon-Perron operator.  Let $F$ be a simple function $\mmu$-almost everywhere, that is, there exist a finite set of real symmetric matrices $\{B_1, \dots, B_m\}$ such that $F^{-1}(B_j)=\{x\in \mathcal{X}\mid F(x)=B_j\}$ satisfies
    \[\mmu(F^{-1}(B_j))>0, \forall j\in [m], \quad \sum_{j\in [m]} \mmu(F^{-1}(B_j))=1.\]
    Assume that for all $B_j$, $aI\preceq B_j\preceq bI$, and let these be tight.  Let $T$ be the operator defined as $T(x)=\frac{1}{2}(F(x)\otimes I_{d}+I_{d}\otimes F(x))$, with exponential multiplication operator $E_T^{\theta}$.  Define the matrix-valued function
    \[\mathcal{F}(r)=\mathbb{E}_{\mmu}\bbb{(1-\lambda)\exp\bb{\frac{\theta}{2}T(x)}\bb{rI-\lambda \exp(\theta T(x))}^{-1}\exp\bb{\frac{\theta}{2}T(x)}     }.\]
    Then the following hold:
    \begin{enumerate}
        \item Let $r^*$ be such that $\mathcal{F}(r^*)$ has an eigenvalue of 1, and let $\mathbf{v}$ be an eigenvector corresponding to the eigenvalue 1 of $\mathcal{F}(r^*)$.  Then $r^*$ is an eigenvalue of $E_T^{\theta/2}\hat{P}E_T^{\theta/2}$ with corresponding eigenfunction $\mathbf{h}$ such that
        \[\mathbf{h}(x)=(1-\lambda)(r^* I-\lambda \exp(\theta T(x)))^{-1}\exp\bb{\frac{\theta}{2} T(x)}\mathbf{v}.\]
        Conversely, every eigenvalue of $E_T^{\theta/2}\hat{P}E_T^{\theta/2}$ is defined this way.
        \item Let $\rho$ be the largest eigenvalue of $E_T^{\theta/2}\hat{P}E_T^{\theta/2}$.  Then $\rho> \lambda e^{\theta b}$.
        \item $\norm{E_T^{\theta/2}\hat{P}E_T^{\theta/2}}_{\mmu}=\rho$ and 
        \[\lim_{n\to\infty}\frac{1}{n}\log \inner{\mathbf{1}\otimes \operatorname{vec}(I_d)}{E^{\theta/2}_T(E^{\theta/2}_T\hat{P}E_T^{\theta/2})^{n-1}E^{\theta/2}_T (\mathbf{1}\otimes \operatorname{vec}(I_d))}_{\mmu}=\log \rho = \log \norm{E_T^{\theta/2}\hat{P}E_T^{\theta/2}}_{\mmu}.\]
    \end{enumerate}
\end{lemma}
\begin{proof}
For (1), if $r^*$ is an eigenvalue of $E_T^{\theta/2}\hat{P}E_T^{\theta/2}$ with eigenfunction $\mathbf{h}$, then
    \begin{align*}
         & E_T^{\theta/2}\hat{P} E_T^{\theta/2}\mathbf{h}-r^* \mathbf{h}=0 \\
         \iff \, & \lambda E^{\theta}_T \mathbf{h} +(1-\lambda)E_T^{\theta/2}\tilde{\Pi}E_T^{\theta/2}\mathbf{h}-r^* \mathbf{h}=0 \\
         \iff \, & \lambda \exp(\theta T(x))\mathbf{h}(x)+(1-\lambda)\exp\bb{\frac{\theta}{2}T(x)}\mathbb{E}_{\mmu}\bbb{\exp\bb{\frac{\theta}{2}T(y)}\mathbf{h}(y)}-r^*\mathbf{h}(x)=0.
    \end{align*}
    Let $\mathbf{v}$ as defined in the statement of part (1).  Plugging in $\mathbf{h}$, the first term on the left-hand side is equal to
    \[\lambda (1-\lambda)\exp(\theta T(x))(r^* I-\lambda \exp(\theta T(x)))^{-1}\exp\bb{\frac{\theta}{2} T(x)}\mathbf{v}.\]
    The second term is equal to 
    \begin{align*}
        & (1-\lambda)\exp\bb{\frac{\theta}{2}T(x)}\mathbb{E}_{\mmu}\bbb{(1-\lambda)\exp\bb{\frac{\theta}{2}T(y)}\bb{r^*I-\lambda \exp(\theta T(y))}^{-1}\exp\bb{\frac{\theta}{2}T(y)}\mathbf{v}}\\
        =\, & (1-\lambda)\exp\bb{\frac{\theta}{2}T(x)}\mathcal{F}(r^*)\mathbf{v}\\
        =\, & (1-\lambda)\exp\bb{\frac{\theta}{2}T(x)}\mathbf{v}\\
        =\, & (1-\lambda)(r^* I-\lambda \exp(\theta T(x)))(r^* I-\lambda \exp(\theta T(x)))^{-1}\exp\bb{\frac{\theta}{2} T(x)}\mathbf{v}\\
        =\, & r^* (1-\lambda)(r^* I-\lambda \exp(\theta T(x)))^{-1}\exp\bb{\frac{\theta}{2} T(x)}\mathbf{v}\\
        &\quad -\lambda(1-\lambda)\exp(\theta T(x)) (r^* I-\lambda \exp(\theta T(x)))^{-1}\exp\bb{\frac{\theta}{2} T(x)}\mathbf{v}
    \end{align*}
    The third term is equal to
    \[r^*(1-\lambda)(r^* I-\lambda \exp(\theta T(x)))^{-1}\exp\bb{\frac{\theta}{2} T(x)}\mathbf{v}.\]
    Adding the first two and subtracting the third gives zero, and shows that $r^*$ is indeed an eigenvalue with eigenfunction $\mathbf{h}$.

    Conversely, if $(r^*, \mathbf{h})$ is an eigenpair of $E_T^{\theta/2}\hat{P}E_T^{\theta/2}$, then solving for $\mathbf{h}$ in the equation 
    \[\lambda \exp(\theta T(x))\mathbf{h}(x)+(1-\lambda)\exp\bb{\frac{\theta}{2}T(x)}\mathbb{E}_{\mmu}\bbb{\exp\bb{\frac{\theta}{2}T(y)}\mathbf{h}(y)}-r^*\mathbf{h}(x)=0\]
    gives
    \[\mathbf{h}(x)=(1-\lambda)(r^* I-\lambda \exp(\theta T(x)))^{-1}\exp\bb{\frac{\theta}{2}T(x)}\mathbb{E}_{\mmu}\bbb{\exp\bb{\frac{\theta}{2}T(x)}\mathbf{h}(x)}.\]
    Multiplying by $\exp\bb{\frac{\theta}{2}T(x)}$ on both sides and taking expectations, we see that
    \[\mathbb{E}_{\mmu}\bbb{\exp\bb{\frac{\theta}{2}T(x)}\mathbf{h}(x)}=\mathcal{F}(r^*)\mathbb{E}_{\mmu}\bbb{\exp\bb{\frac{\theta}{2}T(x)}\mathbf{h}(x)}\]
    so indeed $r^*$ is such that $\mathcal{F}(r^*)$ has an eigenvalue of 1 with $\mathbf{v}=\mathbb{E}_{\mmu}\bbb{\exp\bb{\frac{\theta}{2}T(x)}\mathbf{h}(x)}$.
    
For (2), the function $\mathcal{F}(r)$ is continuous in $r$ when $r$ is large enough (so that the inverses exist); indeed, since $F$ is simple the expectation decomposes as a finite sum, so $\mathcal{F}$ is continuous if all of the summands are continuous.  In particular, since the eigenvalues of all $B_j$ are bounded above by $b$, when $r>\lambda e^{\theta b}$ $\mathcal{F}(r)$ is continuous in $r$.  Now when $r\to \infty$ the entries of $\mathcal{F}(r)$ approach zero, implying that all eigenvalues of $\mathcal{F}(r)$ also approach zero as the eigenvalues are continuous in the entries.  
On the other hand, $\mathcal{F}(r)$ decomposes as a finite sum.  Consider the summand corresponding to some $B_j$ such that $b$ is its largest eigenvalue.  As $r\to\lambda e^{\theta b}$ from above it is clear that the largest eigenvalue of this summand goes to $\infty$; then it is clear that the largest eigenvalue of this finite sum must also go to $\infty$; implying that at least one eigenvalue of $\mathcal{F}(r)$ goes to $\infty$.
Thus, as $\mathcal{F}(r)$ is entrywise continuous in $r$, and as the eigenvalues are continuous in the entries, for some $r^* \in (\lambda e^{\theta b}, \infty)$ $\mathcal{F}(r^*)$ has an eigenvalue of 1.  By part (1). this implies that $r^*$ is an eigenvalue of $E_T^{\theta/2}\hat{P}E_T^{\theta/2}$.  Now let $\rho$ be the largest of all such eigenvalues, and so we have shown that $\rho> \lambda e^{\theta b}$.
    
For (3), let $\sigma(E_T^{\theta/2}\hat{P}E_T^{\theta/2})$ be the spectrum of $E_T^{\theta/2}\hat{P}E_T^{\theta/2}$.  As $E_T^{\theta/2}\hat{P}E_T^{\theta/2}$ is self-adjoint on $\ell_2(\mmu\otimes \mathbf{1})$, the spectrum is real and decomposes as the disjoint union of the essential spectrum $\sigma_{\text{ess}}(E_T^{\theta/2}\hat{P}E_T^{\theta/2})$ and the discrete spectrum $\sigma_\text{d}(E_T^{\theta/2}\hat{P}E_T^{\theta/2})$, the latter consisting of all eigenvalues with finite multiplicity.  Define the essential spectral radius as the largest element of the essential spectrum and the spectral radius as the largest element of the spectrum, so the latter is always at least the former.  Now $E_T^{\theta/2}\hat{P}E_T^{\theta/2}=\lambda E_T^{\theta}+(1-\lambda)E^{\theta/2}\tilde{\Pi}E^{\theta/2}$ and thus is a compact perturbation of the operator $\lambda E^{\theta}_T$, as $(1-\lambda)E^{\theta/2}\tilde{\Pi}E^{\theta/2}$ is finite rank (rank at most $d^2$).  Weyl's perturbation theorem \cite{weyl1909} ensures that $\sigma_{\text{ess}}(E_T^{\theta/2}\hat{P}E_T^{\theta/2})=\sigma_{\text{ess}}(\lambda E^{\theta/2}_T)$.  As $\lambda E_T^{\theta/2}$ is a block diagonal operator (Definition \ref{def:mul_operator}), its spectrum is the union of the spectra of $\lambda \exp(\theta T(x)), x\in \mathcal{X}$ \cite{nauimark1964}, and so the essential spectral radius is bounded above by $\lambda e^{\theta b}$.  Therefore, by part (2). the largest eigenvalue $\rho$ is equal to the spectral radius of $E_T^{\theta/2}\hat{P}E_T^{\theta/2}$ and thus $\norm{E_T^{\theta/2}\hat{P}E_T^{\theta/2}}_{\mmu}=\rho$.

To prove the second statement, from Lemma \ref{lem:P_hat_props} we have
\[\inner{\mathbf{1}\otimes \operatorname{vec}(I_d)}{E_T^{\theta/2}\bb{E^{\theta/2}_T\hat{P}E_T^{\theta/2}}^{n-1} E_T^{\theta/2}(\mathbf{1}\otimes \operatorname{vec}(I_d))}_{\mmu}\leq d\norm{E_T^{\theta/2}\hat{P}E_{T}^{\theta/2}}_{\mmu}^n=d\rho^n.\]
This implies
\begin{equation}\label{eq:lim_sup}
    \limsup_{n\to\infty} \frac{1}{n}\log \inner{\mathbf{1}\otimes \operatorname{vec}(I_d)}{E_T^{\theta/2}\bb{E^{\theta/2}_T\hat{P}E_T^{\theta/2}}^{n-1} E_T^{\theta/2}(\mathbf{1}\otimes \operatorname{vec}(I_d))}_{\mmu}\leq \log \rho.
\end{equation}
        
        Let $\mathbf{h}$ be a unit norm eigenfunction corresponding to $\rho$, and let $\tilde{\mathbf{h}}$ be the projection of $E_T^{\theta/2}(\mathbf{1}\otimes \operatorname{vec}(I_d))$ onto $\mathbf{h}$; in other words,
\begin{equation}\label{eq:h_tilde}
    \tilde{\mathbf{h}}=\inner{E_T^{\theta/2}(\mathbf{1}\otimes \operatorname{vec}(I_d))}{\mathbf{h}}_{\mmu}\mathbf{h}.
\end{equation}
Then
\[\inner{E^{\theta/2}_T(\mathbf{1}\otimes \operatorname{vec}(I_d))-\tilde{\mathbf{h}}}{(E_T^{\theta/2}\hat{P}E_T^{\theta/2})^{n-1}\tilde{\mathbf{h}}}_{\mmu}=0\]
by definition of $\tilde{\mathbf{h}}$ being an eigenfunction and it being a projection of $E_T^{\theta/2}(\mathbf{1}\otimes \operatorname{vec}(I_d))$.  Next, we have
\[\inner{E^{\theta/2}_T(\mathbf{1}\otimes \operatorname{vec}(I_d))-\tilde{\mathbf{h}}}{(E_T^{\theta/2}\hat{P}E_T^{\theta/2})^{n-1}(E_T^{\theta/2}(\mathbf{1}\otimes \operatorname{vec}(I_d))-\tilde{\mathbf{h}})}_{\mmu}\geq 0\]
as $E_T^{\theta/2}\hat{P}E_T^{\theta/2}$ is positive semidefinite on $\ell_2(\mmu\otimes \mathbf{1})$ -- this follows from $E_T^{\theta/2}$ being self-adjoint and $\hat{P}$ being positive semidefinite on $\ell_2(\mmu\otimes \mathbf{1})$.  Then
\begin{align*}
    & \inner{E^{\theta/2}_T(\mathbf{1}\otimes \operatorname{vec}(I_d))}{(E_T^{\theta/2}\hat{P}E_T^{\theta/2})^{n-1}E^{\theta/2}_T(\mathbf{1}\otimes \operatorname{vec}(I_d))}_{\mmu}\\
    =\, & \inner{\tilde{\mathbf{h}}}{(E_T^{\theta/2}\hat{P}E_T^{\theta/2})^{n-1}\tilde{\mathbf{h}}}_{\mmu}\\
    &\quad + 2\inner{E^{\theta/2}_T(\mathbf{1}\otimes \operatorname{vec}(I_d))-\tilde{\mathbf{h}}}{(E_T^{\theta/2}\hat{P}E_T^{\theta/2})^{n-1}\tilde{\mathbf{h}}}_{\mmu}\\
    &\quad + \inner{E^{\theta/2}_T(\mathbf{1}\otimes \operatorname{vec}(I_d))-\tilde{\mathbf{h}}}{(E_T^{\theta/2}\hat{P}E_T^{\theta/2})^{n-1}(E_T^{\theta/2}(\mathbf{1}\otimes \operatorname{vec}(I_d))-\tilde{\mathbf{h}})}_{\mmu}\\
    &\geq \inner{\tilde{\mathbf{h}}}{(E_T^{\theta/2}\hat{P}E_T^{\theta/2})^{n-1}\tilde{\mathbf{h}}}_{\mmu}\\
    &=\inner{E_T^{\theta/2}(\mathbf{1}\otimes \operatorname{vec}(I_d))}{\mathbf{h}}_{\mmu}^2\rho^{n-1}.
\end{align*}
Assume $\inner{E_T^{\theta/2}(\mathbf{1}\otimes \operatorname{vec}(I_d))}{\mathbf{h}}_{\mmu}\neq 0$ for now.  Then the above shows that
\begin{equation}\label{eq:lim_inf}
    \liminf_{n\to\infty} \frac{1}{n}\log \inner{\mathbf{1}\otimes \operatorname{vec}(I_d)}{E_T^{\theta/2}\bb{E^{\theta/2}_T\hat{P}E_T^{\theta/2}}^{n-1} E_T^{\theta/2}(\mathbf{1}\otimes \operatorname{vec}(I_d))}_{\mmu}\geq \log\rho.
\end{equation}
Putting Equations \ref{eq:lim_sup} and \ref{eq:lim_inf} together proves the statement. 

We now show that $\inner{E_T^{\theta/2}(\mathbf{1}\otimes \operatorname{vec}(I_d))}{\mathbf{h}}_{\mmu}\neq 0$.  To do this we make a connection to dynamical systems and the theory of positive operators in Banach spaces, and which uses our realness assumption on $F$.  For a cone $K$ in a Banach space, an operator $A$ is \textit{positive} with respect to $K$ if $AK\subseteq K$.  In our case, let the Banach space under consideration be the space of $d\times d$ symmetric matrices equipped with the Frobenius norm, or equivalently, the space of length $d^2$ vectors whose reshaping to a $d\times d$ matrix is symmetric equipped with the Euclidean norm; we will show that it is no loss to restrict to this space.  Let $\mathbf{S}^d_+$ be the cone of real positive semidefinite $d\times d$ matrices; it is a closed, convex cone.  It is also \textit{reproducing} in the space of real symmetric $d\times d$ matrices, in that every symmetric matrix $X$ can be written as $X=X_1-X_2$ for $X_1, X_2$ positive semidefinite \cite{hill1987}.  We can also identify $\mathbf{S}_+^d$ as the space of length $d^2$ vectors whose reshaping to a $d\times d$ matrix is positive semidefinite; we now make this identification.  Now let $K\subset \ell_2(\mmu\otimes \mathbf{1})$ be the infinite direct product of $\mathbf{S}_+^d$ indexed by elements of $\mathcal{X}$; in other words, for $\mathbf{h}\in \ell_2(\mmu\otimes \mathbf{1})$, $\mathbf{h}\in K$ if and only if $\mathbf{h}(x)\in \mathbf{S}_+^d$ for all $x\in \mathcal{X}$.  This cone is is closed, convex, and reproducing in the space of all $\mathbf{h}\in\ell_2(\mmu\otimes \mathbf{1})$ such that $\mathbf{h}(x)$ can be reshaped to a symmetric $d\times d$ matrix for all $x\in\mathcal{X}$.

Next, we show that we can indeed restrict ourselves to the Banach space of symmetric matrices equipped with the Frobenius norm, in the sense that any eigenfunction $\mathbf{h}$ of $E_T^{\theta/2}\hat{P}E_T^{\theta/2}$ can be such that $\mathbf{h}(x)$ can be reshaped into a $d\times d$ symmetric matrix.  As an operator, we see that
\[(E_T^{\theta/2}\hat{P}E_T^{\theta/2}\mathbf{h})(x)=\int \exp\bb{\frac{\theta}{2}T(x)}\exp\bb{\frac{\theta}{2}T(y)}\mathbf{h}(y) P(x, dy).\]
From the definition of $T$ we see that $\exp\bb{\frac{\theta}{2}T(x)}=\exp\bb{\frac{\theta}{4}F(x)}\otimes \exp\bb{\frac{\theta}{4}F(x)}$, and so $\exp\bb{\frac{\theta}{2}T(x)}\exp\bb{\frac{\theta}{2}T(y)}=\exp\bb{\frac{\theta}{4}F(x)}\exp\bb{\frac{\theta}{4}F(y)}\otimes \exp\bb{\frac{\theta}{4}F(x)}\exp\bb{\frac{\theta}{4}F(y)}$.  If $\mathbf{h}(y)$ has a decomposition $\sum_{j} \mathbf{a}_j(y)\otimes \mathbf{b}_j(y)$ then 
\begin{align*}
    &\exp\bb{\frac{\theta}{2}T(x)}\exp\bb{\frac{\theta}{2}T(y)}\mathbf{h}(y)\\
    =&\sum_j \exp\bb{\frac{\theta}{4}F(x)}\exp\bb{\frac{\theta}{4}F(y)}\mathbf{a}_j(y)\otimes \exp\bb{\frac{\theta}{4}F(x)}\exp\bb{\frac{\theta}{4}F(y)}\mathbf{b}_j(y).
\end{align*}

Then it is not hard to see that the function $\mathbf{h}'$ such that $\mathbf{h}'(y)=\sum_j \mathbf{b}_j(y)\otimes \mathbf{a}_j(y)$ is also an eigenfunction of $E_T^{\theta/2}\hat{P}E_T^{\theta/2}$.  So take $(\mathbf{h}+\mathbf{h}')/2$ as the eigenfunction we desire.  Next, we see that $E_T^{\theta/2}\hat{P}E_T^{\theta/2}$ is positive with respect to $K$.  For $\mathbf{h}\in\ell_2(\mmu\otimes \mathbf{1})$ denote $\hat{\mathbf{h}}(x)$ as the $d\times d$ matrix reshaped from the length $d^2$ vector $\mathbf{h}(x)$.  From the action of $E_T^{\theta/2}\hat{P}E_T^{\theta/2}$ and properties of the Kronecker product we have
\[\exp\bb{\frac{\theta}{2}T(x)}\exp\bb{\frac{\theta}{2}T(y)}\mathbf{h}(y)=\exp\bb{\frac{\theta}{4}F(x)}\exp\bb{\frac{\theta}{4}F(y)}\hat{\mathbf{h}}(y)\exp\bb{\frac{\theta}{4}F(y)}\exp\bb{\frac{\theta}{4}F(x)}.\]
It is clear that since $\hat{\mathbf{h}}(y)$ is positive semidefinite so is the right-hand side above.  Then by linearity we see that the integral of this with respect to $P(x, dy)$ is also positive semidefinite.  Therefore, $E_T^{\theta/2}\hat{P}E_T^{\theta/2}\mathbf{h}\in K$, since $x$ was arbitrary.

Now we use Corollary 2.2 of \cite{nussbaum1981}, which ensures that if a cone is closed, convex, and reproducing, and if $A$ is a positive linear operator with respect to the cone with spectral radius strictly larger than the essential spectral radius, then the spectral radius is an eigenvalue and there exists a corresponding nonzero eigenfunction \textit{that lies in the cone}.  In our case, we indeed have all the above are satisfied: our $K$ is closed, convex, and reproducing, $E_T^{\theta/2}\hat{P}E_T^{\theta/2}$ is positive with respect to $K$, and the spectral radius $\rho$ is strictly larger than the essential spectral radius, which is at most $\lambda e^{\theta b}$.  Then there is a nonzero eigenfunction $\mathbf{h}$ in $K\subseteq \ell_2(\mmu\otimes \mathbf{1})$.

This allows us to assert that in the derivation of Equation \ref{eq:lim_inf} we could have chosen $\mathbf{h}\in K$.  With this choice of $\mathbf{h}$ it is straightforward in showing $\inner{E_T^{\theta/2}(\mathbf{1}\otimes \operatorname{vec}(I_d))}{\mathbf{h}}_{\mmu}\neq 0$.  Indeed, we have
\begin{align*}
    \inner{E_T^{\theta/2}(\mathbf{1}\otimes \operatorname{vec}(I_d))}{\mathbf{h}}_{\mmu}
    &=\inner{\mathbf{1}\otimes \operatorname{vec}(I_d)}{E_T^{\theta/2}\mathbf{h}}_{\mmu}\\
    &=\inner{\operatorname{vec}(I_d)}{\mathbb{E}_{\mmu}\bbb{\exp\bb{\frac{\theta}{2}T(x)}\mathbf{h}(x)}}\\
    &=\operatorname{tr}\mathbb{E}_{\mmu}\bbb{\exp\bb{\frac{\theta}{4}F(x)}\hat{\mathbf{h}}(x)\exp\bb{\frac{\theta}{4}F(x)}}\\
    &=\mathbb{E}_{\mmu}\bbb{\operatorname{tr}\bbb{\exp\bb{\frac{\theta}{4}F(x)}\hat{\mathbf{h}}(x)\exp\bb{\frac{\theta}{4}F(x)}}}.
\end{align*}
Each of these matrices $\exp\bb{\frac{\theta}{4}F(x)}\hat{\mathbf{h}}(x)\exp\bb{\frac{\theta}{4}F(x)}$ are positive semidefinite; moreover, as $\mathbf{h}$ is nonzero in the $\ell_2$ space $\hat{\mathbf{h}}(x)$ is nonzero on a set of positive measure (i.e., $\hat{\mathbf{h}}(x)$ cannot be nonzero only on a set of measure zero).  Therefore, since $\exp\bb{\frac{\theta}{4}F(x)}$ is full rank, by Sylvester's inertia theorem \cite{sylvester1852} we have that $\exp\bb{\frac{\theta}{4}F(x)}\hat{\mathbf{h}}(x)\exp\bb{\frac{\theta}{4}F(x)}$ has at least one strictly positive eigenvalue for all $x$ in a set of positive measure, meaning that the trace is strictly positive.  Then \[\inner{E_T^{\theta/2}(\mathbf{1}\otimes \operatorname{vec}(I_d))}{\mathbf{h}}_{\mmu}>0.\]
    
\end{proof}

Using this lemma, we are now able to give the proof of Lemma \ref{lem:limit}.
\begin{proof}[Proof of Lemma \ref{lem:limit}]
    We first reduce to the simple case.  The range of $F$ is compact, as $F$ can be seen as a continuous map of the compact set $O(d)\times [a, b]$.  Then for $\epsilon>0$ let $\mathcal{N}_\epsilon$ be an $\epsilon$-net of the range of $F$ with respect to the metric induced by the operator norm.  Define $F_\epsilon$ as $F_\epsilon(x)=\argmin_{B\in \mathcal{N}_\epsilon}\norm{F(x)-B}$.  Then it is clear that for all $x\in \mathcal{X}$, $\norm{F_\epsilon(x)-F(x)}\leq \epsilon$.  Define $T_\epsilon(x)=\frac{1}{2}(F_\epsilon(x)\otimes I_d+I_d\otimes F_\epsilon(x))$ and $E_{T_{\epsilon}}^{\theta/2}$ the corresponding block diagonal operator as usual.  From this we see that $E_{T_\epsilon}^\theta$ converges to $E_T^{\theta}$ in norm, and thus $E_{T_\epsilon}^{\theta/2}\hat{P}E_{T_\epsilon}^{\theta/2}$ converges to $E_{T}^{\theta/2}\hat{P}E_{T}^{\theta/2}$ in norm via any sequence of $\epsilon$-nets with $\epsilon\to 0$.  Then $\displaystyle\lim_{\epsilon\to 0}\log \norm{E_{T_\epsilon}^{\theta/2}\hat{P}E_{T_\epsilon}^{\theta/2}}_{\mmu}=\log \norm{E_{T}^{\theta/2}\hat{P}E_{T}^{\theta/2}}_{\mmu}$.

    From Lemma \ref{lem:simple_function} we have
    \[\lim_{n\to\infty}\frac{1}{n}\log \inner{\mathbf{1}\otimes \operatorname{vec}(I_d)}{E^{\theta/2}_{T_\epsilon}(E^{\theta/2}_{T_\epsilon}\hat{P}E_{T_\epsilon}^{\theta/2})^{n-1}E^{\theta/2}_{T_\epsilon} (\mathbf{1}\otimes \operatorname{vec}(I_d))}_{\mmu}= \log \norm{E_{T_\epsilon}^{\theta/2}\hat{P}E_{T_\epsilon}^{\theta/2}}_{\mmu},\]
    so that, denoting $a_{\epsilon, n}=\frac{1}{n}\log \inner{\mathbf{1}\otimes \operatorname{vec}(I_d)}{E^{\theta/2}_{T_\epsilon}(E^{\theta/2}_{T_\epsilon}\hat{P}E_{T_\epsilon}^{\theta/2})^{n-1}E^{\theta/2}_{T_\epsilon} (\mathbf{1}\otimes \operatorname{vec}(I_d))}_{\mmu}$,
    \[\lim_{\epsilon\to 0}\lim_{n\to\infty} a_{\epsilon, n}= \lim_{\epsilon\to 0}\log \norm{E_{T_\epsilon}^{\theta/2}\hat{P}E_{T_\epsilon}^{\theta/2}}_{\mmu}=\log\norm{E_{T}^{\theta/2}\hat{P}E_{T}^{\theta/2}}_{\mmu}.\]
    In fact, we can give a more refined proof of part (3). of Lemma \ref{lem:simple_function} to show that $\displaystyle\lim_{n\to\infty} a_{\epsilon, n}$ uniformly converges in $\epsilon$.  Let $\rho_\epsilon$ be the leading eigenvalue of $E_{T_\epsilon}^{\theta/2}\hat{P} E_{T_\epsilon}^{\theta/2}$.  We first have by Lemma \ref{lem:P_hat_props} that
    \[\inner{\mathbf{1}\otimes \operatorname{vec}(I_d)}{E^{\theta/2}_{T_\epsilon}(E^{\theta/2}_{T_\epsilon}\hat{P}E_{T_\epsilon}^{\theta/2})^{n-1}E^{\theta/2}_{T_\epsilon} (\mathbf{1}\otimes \operatorname{vec}(I_d))}_{\mmu}\leq d\norm{E_{T_\epsilon}^{\theta/2}\hat{P} E_{T_\epsilon}^{\theta/2}}_{\mmu}^n=d\rho_\epsilon^n,\]
    so $a_{\epsilon, n}-\log \rho_\epsilon\leq (\log d)/n$.  Next, we showed for any unit norm eigenfunction $\mathbf{h}_\epsilon$ corresponding to the leading eigenvalue $\rho_\epsilon$ that 
    \[\inner{\mathbf{1}\otimes \operatorname{vec}(I_d)}{E^{\theta/2}_{T_\epsilon}(E^{\theta/2}_{T_\epsilon}\hat{P}E_{T_\epsilon}^{\theta/2})^{n-1}E^{\theta/2}_{T_\epsilon} (\mathbf{1}\otimes \operatorname{vec}(I_d))}_{\mmu}\geq \inner{E_{T_\epsilon}^{\theta/2}(\mathbf{1}\otimes \operatorname{vec}(I_d))}{\mathbf{h}_\epsilon}_{\mmu}^2 \rho_\epsilon^{n-1},\]
    so $a_{\epsilon, m}-\log \rho_\epsilon\geq \bb{\log \inner{E_{T_\epsilon}^{\theta/2}(\mathbf{1}\otimes \operatorname{vec}(I_d))}{\mathbf{h}}_{\mmu}^2}/n-(\log \rho_\epsilon)/n$. 

    In the proof of part (3). of Lemma \ref{lem:simple_function} we showed that $\mathbf{h}_\epsilon$ can in fact be chosen such that for all $x\in\mathcal{X}$, $\mathbf{h}_\epsilon(x)$ can be reshaped to a positive semidefinite $d\times d$ matrix.  With such structure we lower bounded $\inner{E_{T_\epsilon}^{\theta/2}(\mathbf{1}\otimes \operatorname{vec}(I_d))}{\mathbf{h}_\epsilon}_{\mmu}>0$.  We can actually give a more quantitative lower bound of this quantity that is independent of $\epsilon$; in addition, we can give an upper bound $\rho_\epsilon$ also independent of $\epsilon$.  First note that under an $\epsilon$-net of the range of $F$, since $aI\preceq F(x)\preceq bI$, that $aI\preceq F_{\epsilon}(x)\preceq bI$ also.  Then the upper bound is simple; we have $\rho_\epsilon=\norm{E_{T_{\epsilon}}^{\theta/2}\hat{P}E_{T_\epsilon}^{\theta/2}}_{\mmu}\leq e^{\theta b}$.  For the lower bound we first write
    \begin{align*}
        \inner{E_{T_\epsilon}^{\theta/2}(\mathbf{1}\otimes \operatorname{vec}(I_d))}{\mathbf{h}_\epsilon}_{\mmu}
        &=\inner{\operatorname{vec}(I_d)}{\mathbb{E}_{\mmu}\bbb{\exp\bb{\frac{\theta}{2}T_{\epsilon}(x)}\mathbf{h}_\epsilon(x)}}\\
        &=\operatorname{tr}\mathbb{E}_{\mmu}\bbb{\exp\bb{\frac{\theta}{4}F_\epsilon(x)}\hat{\mathbf{h}_\epsilon}(x)\exp\bb{\frac{\theta}{4}F_\epsilon(x)}}.
    \end{align*}
    Since we have established $\exp\bb{\frac{\theta}{2}F_\epsilon(x)}\hat{\mathbf{h}_\epsilon}(x)\exp\bb{\frac{\theta}{2}F_\epsilon(x)}$ is positive semidefinite for all $x\in \mathcal{X}$, the trace is equal to the nuclear norm, which upper bounds the Frobenius norm.  Therefore,
    \begin{align*}
        \operatorname{tr}\mathbb{E}_{\mmu}\bbb{\exp\bb{\frac{\theta}{4}F_\epsilon(x)}\hat{\mathbf{h}_\epsilon}(x)\exp\bb{\frac{\theta}{4}F_\epsilon(x)}}
        &= \mathbb{E}_{\mmu}\bbb{\norm{\exp\bb{\frac{\theta}{4}F_\epsilon(x)}\hat{\mathbf{h}_\epsilon}(x)\exp\bb{\frac{\theta}{4}F_\epsilon(x)}}_{*}}\\
        &\geq \mathbb{E}_{\mmu}\bbb{\norm{\exp\bb{\frac{\theta}{4}F_\epsilon(x)}\hat{\mathbf{h}_\epsilon}(x)\exp\bb{\frac{\theta}{4}F_\epsilon(x)}}_{F}}\\
        &\geq \mathbb{E}_{\mmu}\bbb{e^{\theta a/2}\norm{\hat{\mathbf{h}_\epsilon}(x)}_{F}}\\
        &=e^{\theta a/2}\mathbb{E}_{\mmu}\bbb{\norm{\mathbf{h}_\epsilon(x)}}.
    \end{align*}
    Now since $F_\epsilon$ takes only finite number of values, we have that $\mathbf{h}_\epsilon$ can also only take a finite number of values (see part (1) of Lemma \ref{lem:simple_function}).  Letting $A_j=F^{-1}(B_j)$ and denoting $\mathbf{h}_j=\mathbf{h}_\epsilon(x)$ for $x\in A_j$, we have
    \[\mathbb{E}_{\mmu}\bbb{\norm{\mathbf{h}_\epsilon(x)}}=\sum_{j\in [m]}\mmu(A_j)\norm{\mathbf{h}_j}.\]
    Now as $\norm{\mathbf{h}_\epsilon}_{\mmu}=1$ we must have $\norm{\mathbf{h}_j}\leq 1$, and so then $\norm{\mathbf{h}_j}\geq \norm{\mathbf{h}_j}^2$.  But $\sum_{j\in [m]}\mmu(A_j)\norm{\mathbf{h}_j}^2=\norm{\mathbf{h}_\epsilon}_{\mmu}^2$, and thus $\mathbb{E}_{\mmu}[\norm{\mathbf{h}_\epsilon(x)}]\geq 1$.

    Therefore, $\inner{E_{T_\epsilon}^{\theta/2}(\mathbf{1}\otimes \operatorname{vec}(I_d))}{\mathbf{h}}_{\mmu}^2\geq e^{\theta a}$, and we thus have $a_{\epsilon, n}-\log \rho_\epsilon\geq (\theta(a-b))/n$.  Putting this all together, we have shown
    \[\abs{a_{\epsilon, n}-\log \rho_\epsilon}\leq \max\left\{\frac{\log d}{n}, \frac{\theta(b-a)}{n}\right\},\]
    which implies uniform convergence, since $d, \theta, a$, and $b$ are all constants independent of $\epsilon$.  We can then swap limits:
    \[\lim_{\epsilon\to 0}\lim_{n\to\infty} a_{\epsilon, n}=\lim_{n\to\infty}\lim_{\epsilon\to 0}a_{\epsilon, n}.\]
    But as $E_{T_{\epsilon}}^\theta$ converges to $E_T^{\theta/2}$ in norm, so does $E^{\theta/2}_{T_\epsilon}(E^{\theta/2}_{T_\epsilon}\hat{P}E_{T_\epsilon}^{\theta/2})^{n-1}E^{\theta/2}_{T_\epsilon}$ to $E^{\theta/2}_{T}(E^{\theta/2}_{T}\hat{P}E_{T}^{\theta/2})^{n-1}E^{\theta/2}_{T}$, which implies weak convergence, that is, 
    \begin{align*}
        & \lim_{\epsilon \to 0}\inner{\mathbf{1}\otimes \operatorname{vec}(I_d)}{E^{\theta/2}_{T_\epsilon}(E^{\theta/2}_{T_\epsilon}\hat{P}E_{T_\epsilon}^{\theta/2})^{n-1}E^{\theta/2}_{T_\epsilon} (\mathbf{1}\otimes \operatorname{vec}(I_d))}_{\mmu}\\
        =\, & \inner{\mathbf{1}\otimes \operatorname{vec}(I_d)}{E^{\theta/2}_{T}(E^{\theta/2}_{T}\hat{P}E_{T}^{\theta/2})^{n-1}E^{\theta/2}_{T} (\mathbf{1}\otimes \operatorname{vec}(I_d))}_{\mmu}.
    \end{align*}
    Then we have the chain of equalities
    \begin{align*}
        & \lim_{n\to\infty}\frac{1}{n}\log \mathbb{E}_{\mmu}\bbb{\norm{\prod_{j=1}^n \exp\bb{\frac{\theta }{2}F(s_j)}}_F^2}\\
        =\, & \lim_{n\to\infty}\frac{1}{n}\log \inner{\mathbf{1}\otimes \operatorname{vec}(I_d)}{E^{\theta/2}_{T}(E^{\theta/2}_{T}\hat{P}E_{T}^{\theta/2})^{n-1}E^{\theta/2}_{T} (\mathbf{1}\otimes \operatorname{vec}(I_d))}_{\mmu}\\
        =\, & \lim_{n\to\infty}\lim_{\epsilon\to 0}\frac{1}{n}\log \inner{\mathbf{1}\otimes \operatorname{vec}(I_d)}{E^{\theta/2}_{T_\epsilon}(E^{\theta/2}_{T_\epsilon}\hat{P}E_{T_\epsilon}^{\theta/2})^{n-1}E^{\theta/2}_{T_\epsilon} (\mathbf{1}\otimes \operatorname{vec}(I_d))}_{\mmu}\\
        =\, & \lim_{\epsilon\to 0}\lim_{n\to\infty}\frac{1}{n}\log \inner{\mathbf{1}\otimes \operatorname{vec}(I_d)}{E^{\theta/2}_{T_\epsilon}(E^{\theta/2}_{T_\epsilon}\hat{P}E_{T_\epsilon}^{\theta/2})^{n-1}E^{\theta/2}_{T_\epsilon} (\mathbf{1}\otimes \operatorname{vec}(I_d))}_{\mmu}\\
        =\, & \lim_{\epsilon\to 0} \log \norm{E_{T_\epsilon}^{\theta/2}\hat{P} E_{T_\epsilon}^{\theta/2}}_{\mmu}\\
        =\, & \log \norm{E_{T}^{\theta/2}\hat{P} E_{T}^{\theta/2}}_{\mmu}.
    \end{align*}
\end{proof}

\section{Final bounds}\label{sec:final}
In this section we prove our main theorems regarding the moment generating function and give concentration inequalities.  We first show how to go from bounds on the moment generating function of Equation \ref{eq:matrix_mgf} to tail bounds.  

\subsection{From the  MGF to tail bounds}
The reason we have studied the moment generating function of Equation \ref{eq:matrix_mgf} is due to the following result:
\begin{theorem}[Multi-matrix Golden-Thompson inequality, \cite{garg2018}]\label{thm:mmgm}
Let $H_1, \dots, H_k\in \mathbb{C}^{d\times d}$ be Hermitian matrices.  Then
\[\log\bb{\operatorname{tr}\bbb{\exp\bb{\sum_{j=1}^k H_j}}}\leq \frac{4}{\pi}\int_{-\frac{\pi}{2}}^{\frac{\pi}{2}}\log \bb{\operatorname{tr}\bbb{\prod_{j=1}^k \exp\bb{\frac{e^{i\phi}}{2}H_j}\prod_{j=k}^1\exp\bb{\frac{e^{-i\phi}}{2}H_j}}}\,d\boldsymbol{\xi}(\phi),\]
where $\boldsymbol{\xi}$ is some probability distribution on $\bbb{-\frac{\pi}{2}, \frac{\pi}{2}}$.
\end{theorem}

From this, we can obtain tail bounds.  Recall that $F_j:\mathcal{X}\to \mathbb{R}^{d\times d}$ are a sequence of measurable functions from $\mathcal{X}$ to $d\times d$ real symmetric matrices, and let $\mathbb{E}[F_j]=0$.  First, 
\begin{equation}
    \begin{aligned}
    \operatorname{Pr}_{\mmu}\bb{\lambda_{\text{max}}\bb{\sum_j F_j(s_j)}\geq t}
    &\leq e^{-\theta t}\mathbb{E}_{\mmu}\bbb{\exp\bb{\theta \lambda_{\text{max}}\bb{\sum_j F_j(s_j)}}}\\
    &= e^{-\theta t}\mathbb{E}_{\mmu}\bbb{\lambda_{\text{max}}\bb{\exp\bb{\theta \sum_j F_j(s_j)}}}\\
    &\leq e^{-\theta t}\mathbb{E}_{\mmu}\bbb{\operatorname{tr}\exp\bb{\theta \sum_j F_j(s_j)}}.
    \end{aligned}
\end{equation}
We would like to apply Theorem \ref{thm:mmgm}; we follow the standard outline given in \cite{garg2018}.  An immediate application of Jensen's inequality on the right-hand side of Theorem \ref{thm:mmgm} furnishes an upper bound, and then taking exponents, we have
\[\operatorname{tr}\bbb{\exp\bb{\sum_{j=1}^n F_j(s_j)}}\leq \bb{\int_{-\frac{\pi}{2}}^{\frac{\pi}{2}}\operatorname{tr}\bbb{\prod_{j=1}^n \exp\bb{\frac{e^{i\phi}}{2}F_j(s_j)}\prod_{j=n}^1\exp\bb{\frac{e^{-i\phi}}{2}F_j(s_j)}}\,d\boldsymbol{\xi}(\phi)}^{4/\pi}.\]
Using $\norm{x}_p\leq d^{1/p-1}\norm{x}_1$ for $p\in (0,1)$ and choosing $p=\pi/4$ we have
\[\operatorname{tr}\bbb{\exp\bb{\frac{\pi}{4}\sum_{j=1}^n F_j(s_j)}}^{4/\pi}\leq d^{4/\pi-1}\operatorname{tr}\bbb{\exp\bb{\sum_{j=1}^n F_j(s_j)}}.\]
Combining the above two lines, adding in $\theta$, and taking expectation, we have
\begin{equation}
\begin{aligned}
    & \mathbb{E}\bbb{\operatorname{tr}\exp\bb{\frac{\pi}{4}\theta \sum_{j=1}^n F_j(s_j)}}\\
    \leq &d^{1-\pi/4}\int_{-\frac{\pi}{2}}^{\frac{\pi}{2}}\mathbb{E}\bbb{\operatorname{tr}\bbb{\prod_{j=1}^n \exp\bb{\frac{\theta e^{i\phi}}{2}F_j(s_j)}\prod_{j=n}^1\exp\bb{\frac{\theta e^{-i\phi}}{2}F_j(s_j)}}}\,d\boldsymbol{\xi}(\phi)\\
    =&d^{1-\pi/4}\int_{-\frac{\pi}{2}}^{\frac{\pi}{2}}\mathbb{E}\bbb{\norm{\prod_{j=1}^n \exp\bb{\frac{\theta e^{i\phi}}{2}}F_j(s_j)}_F^2}\,d\boldsymbol{\xi}(\phi).
\end{aligned}
\end{equation}
Now let $M_F(\theta)$ be such that
\begin{equation}\label{eq:M_F(theta)}
    M_F(\theta)\geq \mathbb{E}\bbb{\norm{\prod_{j=1}^n \exp\bb{\frac{\theta e^{i\phi}}{2}}F_j(s_j)}_F^2}
\end{equation}
that is independent of $\phi$, i.e., an upper bound on the moment generating function in Equation \ref{eq:matrix_mgf} valid for any $\phi\in [-\pi/2, \pi/2]$.  Then 
\begin{equation}\label{eq:concentration_bound_general}
    \begin{aligned}
    \operatorname{Pr}_{\mmu}\bb{\lambda_{\text{max}}\bb{\sum_j F_j(s_j)}\geq t}
    &\leq e^{-\theta t}\mathbb{E}_{\mmu}\bbb{\operatorname{tr}\exp\bb{\theta \sum_j F_j(s_j)}}\\
    &\leq d^{1-\pi/4}e^{-\theta t}\int_{\frac{-\pi}{2}}^{\frac{\pi}{2}}M_F\bb{\frac{4}{\pi}\theta}\, d\boldsymbol{\xi}(\phi)\\
    &=d^{1-\pi/4}\exp\bb{-\theta t+\log M_F\bb{\frac{4}{\pi}\theta}},
    \end{aligned}
\end{equation}
where the last line holds because $\boldsymbol{\xi}$ is a probability distribution and $M_F$ does not depend on $\phi$.
 
We now give our final bounds on the moment generating function 
\[\mathbb{E}_{\mmu}\bbb{\norm{\prod_{j=1}^n \exp\bb{\frac{\theta e^{i\phi}}{2}}F_j(s_j)}_F^2}\]
under both Hoeffding and Bernstein-type assumptions.  In Sections \ref{sec:starting} and \ref{sec:operator} we have reduced the problem to estimating the leading eigenvalues of the operators
\begin{equation}\label{eq:main_matrix}
E_{T_j}^{\theta/2}\hat{P}E_{T_j}^{\theta/2}
\end{equation}
for $j=1, \dots, n$ and $\hat{P}$ the lifted Leon-Perron version of $P$.  Lemma \ref{lem:limit} has shown that each eigenvalue is a limit a corresponding moment generating function.  Thus, the general idea will be to bound 
\[\mathbb{E}_{\mmu}\bbb{\norm{\prod_{k=1}^n \exp\bb{\frac{\theta \cos(\phi)}{2} F_j(s_k)}}_F^2}\]
for each $j$, and where $s_1, \dots, s_n$ is a sequence of states driven by the Markov chain corresponding to $\hat{P}$.  In these bounds we will proceed by way of defining a function $M_F$ to bound these operators.

\subsection{Hoeffding}
In this section, fix a matrix-valued function $F$, and let $\mathbb{E}_{\mmu}[F(x)]=0$ and $aI\preceq F(x)\preceq bI$ for all $x\in \mathcal{X}$.  
For a scalar convex functions $\Psi$, we have the following proposition:
\begin{proposition}\label{prop:matrix_convexity}
    Let $H$ be a Hermitian $d\times d$ matrix such that $aI\preceq H\preceq bI$.  Then for any scalar convex function $\Psi$, 
    \[\Psi(H)\preceq \frac{bI-H}{b-a}\Psi(a)+\frac{H-aI}{b-a}\Psi(b).\]
    Here $\Psi$ acts spectrally; if $H=U\Lambda U^*$, then $\Psi(H)=U\Psi(\Lambda)U^*$, where $\Psi$ acts on the diagonal matrix $\Lambda$ diagonally, i.e., $\Psi(\Lambda)$ is the diagonal matrix such that $\Psi(\Lambda)_{ii}=\Psi(\Lambda_{ii})$.
\end{proposition}
\begin{proof}
    We simply need to show that for any unit vector $\mathbf{u}$ that
    \[\mathbf{u}^* \Psi(H)\mathbf{u}\leq \frac{b-\mathbf{u}^* H\mathbf{u}}{b-a}\Psi(a)+\frac{\mathbf{u}^*H\mathbf{u}-aI}{b-a}\Psi(b).\]
    It suffices to show this for eigenvectors $\mathbf{v}$ of $H$.  If $\lambda$ is the corresponding eigenvalue, then $\mathbf{v}^* \Psi(H)\mathbf{v}=\Psi(\lambda)$.  The right-hand side is
    \[\frac{b-\lambda}{b-a}\Psi(a)+\frac{\lambda-a}{b-a}\Psi(b).\]
    Since $\Psi$ is a scalar convex function and $a\leq \lambda\leq b$, 
    \[\Psi(\lambda)\leq \frac{b-\lambda}{b-a}\Psi(a)+\frac{\lambda-a}{b-a}\Psi(b).\]
    For general $\mathbf{u}$, let $\mathbf{u}=\sum_{k=1}^d c_k \mathbf{v}_k$ be the representation of $\mathbf{u}$ in the basis of eigenvectors of $H$.  Then
    \begin{align*}
        \mathbf{u}^* \Psi(H)\mathbf{u}
        &=\sum_{k=1}^d |c_k|^2 \Psi(\lambda_k)\\
        &\leq \sum_{k=1}^d |c_k|^2 \bb{\frac{b-\lambda_k}{b-a}\Psi(a)+\frac{\lambda_k-a}{b-a}\Psi(b)}\\
        &=\frac{b-\mathbf{u}^* H\mathbf{u}}{b-a}\Psi(a)+\frac{\mathbf{u}^* H\mathbf{u}-a}{b-a}\Psi(b).
    \end{align*}
\end{proof}

If in the above $H$ is a random matrix such that $\mathbb{E}[H]=0$, then taking expectations we have
\[\mathbb{E}[H]\preceq \frac{b}{b-a}\Psi(aI)+\frac{-a}{b-a}\Psi(bI),\]
the right-hand side of which defines an implicit two-point distribution on matrices.  Introduce some notation; let $\boldsymbol\mu=[b/(b-a), -a/(b-a)]^\top$ and $Q=\lambda I + (1-\lambda)\mathbf{1}\boldsymbol\mu^\top$.  Then $Q$ is a transition matrix on two states $\{0, 1\}$ with stationary distribution $\boldsymbol\mu$.  Let $G$ be a function such that $G(0)=aI_d$ and $G(1)=bI_d$.

This allows us to make the following argument about a convex ordering of distributions, with proof in the appendix:
\begin{lemma}\label{lem:markov_matrix_convex_majorization}
    Let $F$ be a function from $\mathcal{X}$ to Hermitian $d\times d$ matrices, and assume $\mathbb{E}_{\mmu}[F]=0$ and $aI\preceq F(x)\preceq bI$ for all $x\in \mathcal{X}$.  Let $s_1, \dots, s_n$ be driven by the Leon-Perron operator $\hat{P}$.  Let $Q, \boldsymbol\mu, G$ be as above.  Let $y_1, \dots, y_n$ be driven by $Q$.  Then for any scalar nonnegative convex function $\Psi$, 
    \[\mathbb{E}_{\mmu}\bbb{\norm{\prod_{j=1}^n \Psi(F(s_j))}_F^2}\leq \mathbb{E}_{\boldsymbol\mu}\bbb{\norm{\prod_{j=1}^n \Psi(G(y_j))}_F^2}.\]
\end{lemma}

\begin{remark}
    The condition on $\Psi$ above does not seem to be very strict.  Indeed, it applies to many natural matrix-valued functions, such as even matrix powers and the matrix exponential.
\end{remark}

Applying this to the scalar convex function $x\mapsto \exp\bb{\frac{\theta \cos(\phi)}{2}x}$ yields 
\begin{equation}\label{eq:final_majorization}
    \mathbb{E}_{\mmu}\bbb{\norm{\prod_{j=1}^n \exp\bb{\frac{\theta \cos(\phi)}{2}F(s_j)}}_F^2}\leq \mathbb{E}_{\boldsymbol\mu}\bbb{\norm{\prod_{j=1}^n \exp\bb{\frac{\theta \cos(\phi)}{2}G(y_j)}}_F^2}.
\end{equation}
Now let $M^{\theta}\in\bbR^{2d^2\times 2d^2}$ be the operator on the two-state state space $\{0, 1\}$ such that $(M^\theta \mathbf{h})(0)=\exp(\theta \cos(\phi)a)I_{d^2}$ and $(M^\theta \mathbf{h})(1)=\exp(\theta \cos(\phi)b)I_{d^2}$, i.e., 
\[M^\theta=
\begin{bmatrix}
\exp(\theta \cos(\phi)a)I_{d^2} & 0\\
0 & \exp(\theta \cos(\phi)b)I_{d^2}
\end{bmatrix}.
\]
Then Lemma \ref{lem:limit} indicates that 
\begin{equation}
    \lim_{n\to\infty}\frac{1}{n}\log \mathbb{E}_{\boldsymbol\mu}\bbb{\norm{\prod_{j=1}^n \exp\bb{\frac{\theta \cos(\phi)}{2}G(y_j)}}_F^2}=\log \eta,
\end{equation}
where $\eta$ is the largest eigenvalue of $M^{\theta/2}\hat{Q}M^{\theta/2}$, $\hat{Q}=Q\otimes I_{d^2}$.  Now note that 
\begin{equation}\label{eq:K_theta}
\begin{aligned}
    & M^{\theta/2}\hat{Q}M^{\theta/2}\\
    =\,& \begin{bmatrix}\exp\bb{\frac{\theta \cos(\phi)}{2}a} & 0\\
0 & \exp\bb{\frac{\theta \cos(\phi)}{2}b}\end{bmatrix}Q\begin{bmatrix}\exp\bb{\frac{\theta \cos(\phi)}{2}a} & 0\\
0 & \exp\bb{\frac{\theta \cos(\phi)}{2}b}\end{bmatrix}\otimes I_{d^2}\\
=:\, & K^{\theta}\otimes I_{d^2}. 
\end{aligned}
\end{equation}
Then by properties of the Kronecker product, the eigenvalues of the matrix on the left-hand side are equal to the eigenvalues of $K^{\theta}$, a $2\times 2$ matrix, the largest of which can be solved for explicitly.

Putting this all together,
\begin{equation}
    \begin{aligned}
        \log \rho
        &=\lim_{n\to\infty} \frac{1}{n}\log \mathbb{E}_{\mmu}\bbb{\norm{\prod_{j=1}^n \exp\bb{\frac{\theta \cos(\phi)}{2}F(s_j)}}_F^2}\\
        &\leq \lim_{n\to\infty}\frac{1}{n}\log \mathbb{E}_{\boldsymbol\mu}\bbb{\norm{\prod_{j=1}^n \exp\bb{\frac{\theta \cos(\phi)}{2}G(y_j)}}_F^2}\\
        &=\log \eta\\
        &=\log \lambda_{\text{max}}(K^{\theta}).
    \end{aligned}
\end{equation}

With the above we can establish the following result, with proof in the appendix. 
\begin{lemma}\label{lem:hoeffding_main}
    Let $K^{\theta}$ be the $2\times 2$ matrix defined in Equation \ref{eq:K_theta}, and let $\eta$ be its largest eigenvalue.  Then
    \[\eta\leq \tilde{\eta}(\theta):=\exp\bb{\alpha(\lambda)\cdot \frac{\theta^2\cos^2(\phi)(b-a)^2}{8}}\]
    where $\alpha:\lambda\mapsto (1+\lambda)/(1-\lambda)$.
\end{lemma}

Putting our work above together, we can now prove the two statements of Theorem \ref{thm:markov_matrix_hoeffding}.
\begin{proof}[Proof of Theorem \ref{thm:markov_matrix_hoeffding}]
    The discussion in Lemma \ref{lem:P_hat_props} (in particular Equation \ref{eq:operator_bound}) and Proposition \ref{prop:T_operator} reveal
    \begin{equation}
        \mathbb{E}_{\mmu}\bbb{\norm{\prod_{j=1}^n \exp\bb{\frac{\theta e^{i\phi}}{2}F_j(s_j)}}_F^2}
        \leq d\prod_{j=1}^n \norm{E_{T_j}^{\theta/2}\hat{P}E_{T_j}^{\theta/2}}_{\mmu},
    \end{equation}
    where recall $T_j(x)=\frac{\cos(\phi)}{2}(F_j(x)\otimes I_d+I_d\otimes F_j(x))$, $E_{T_j}^{\theta/2}$ is the multiplication operator with respect to the matrix-valued function $T_j$ and $\hat{P}=(\lambda I + (1-\lambda)\Pi)\otimes I_{d^2}$.  Then Lemma \ref{lem:limit} and Lemma \ref{lem:markov_matrix_convex_majorization} allow us to bound each of these norms with the norms of corresponding matrices arising from a two state chain, each of which can be bounded using Lemma \ref{lem:hoeffding_main}.  In other words, for each $j$ we have, verifying that the operators $T_j$ satisfy the assumptions of the lemmas using Proposition \ref{prop:T_properties}, 
    \begin{equation}
         \norm{E_{T_j}^{\theta/2}\hat{P}E_{T_j}^{\theta/2}}_{\mmu}\leq \exp\bb{\alpha(\lambda)\cdot \frac{\theta^2\cos(\phi)^2(b_j-a_j)^2}{8}}.
    \end{equation}
    Then 
    \begin{equation}
        \mathbb{E}_{\mmu}\bbb{\norm{\prod_{j=1}^n \exp\bb{\frac{\theta e^{i\phi}}{2}F_j(s_j)}}_F^2}\leq d\exp\bb{\theta^2\cdot \alpha(\lambda)\cdot \frac{\cos(\phi)^2\sum_{j=1}^n (b_j-a_j)^2}{8}},
    \end{equation}
    and setting $\phi=0$ we obtain the first statement of the theorem.

    To turn this into a tail bound, since $\cos(\phi)^2\leq 1$, we can let $M_F(\theta)$, as described in Equation \ref{eq:M_F(theta)}, be
    \begin{equation}
        M_F(\theta):=d\exp\bb{\theta^2\cdot \alpha(\lambda)\cdot \frac{\sum_{j=1}^n (b_j-a_j)^2}{8}}.
    \end{equation}
    Using Equation \ref{eq:concentration_bound_general}, we then have
    \begin{equation}
    \begin{aligned}
    \operatorname{Pr}_{\mmu}\bb{\lambda_{\text{max}}\bb{\sum_{j=1}^n F_j(s_j)}\geq t}
    &\leq\inf_{\theta>0} d^{1-\pi/4}\exp\bb{-\theta \cdot t +\log M_F\bb{\frac{4}{\pi}\theta}}\\
    &=\inf_{\theta>0}d^{2-\pi/4}\exp\bb{-\theta \cdot t+\theta^2\cdot \alpha(\lambda)\cdot \frac{2\sum_{j=1}^n (b_j-a_j)^2}{\pi^2}}.
    \end{aligned}
    \end{equation}
    Optimizing this quadratic over $\theta$ we have
    \begin{equation}
        \operatorname{Pr}_{\mmu}\bb{\lambda_{\text{max}}\bb{\sum_{j=1}^n F_j(s_j)}\geq t}
        \leq d^{2-\pi/4}\exp\bb{\frac{-t^2/(8/\pi^2)}{\alpha(\lambda)\cdot \sum_{j=1}^n (b_j-a_j)^2}},
    \end{equation}
    as desired.
\end{proof}

\subsection{Bernstein}
In this section, fix a matrix-valued function $F$, and let $\mathbb{E}_{\mmu}[F(x)]=0$, $\norm{\mathbb{E}_{\mmu}[F(x)^2]}\leq \mathcal{V}$, and $\norm{F(x)}\leq \mathcal{M}$ for all $x\in\mathcal{X}$.

We take a slightly different approach to bounding the norms of the operators $E_{T_j}^{\theta/2}\hat{P}E_{T_j}^{\theta/2}$.  In particular, consider again Equation \ref{eq:matrix_mgf} for a time-independent function $T$ and the Leon-Perron operator $\hat{P}$.  Then
\begin{equation}
    \begin{aligned}
        & \mathbb{E}_{\mmu}\bbb{\norm{\prod_{j=1}^n \exp\bb{\frac{\theta\cos(\phi)}{2}}F(s_j)}_F^2}\\
        =\, & \inner{\mathbf{1}\otimes\operatorname{vec}(I_d)}{E^{\theta/2}_T(E_T^{\theta/2}\hat{P}E_T^{\theta/2})^{n-1}E_T^{\theta/2}(\mathbf{1}\otimes\operatorname{vec}(I_d))}_{\mmu}\\
        =\, & \inner{\tilde{\Pi}(\mathbf{1}\otimes\operatorname{vec}(I_d))}{E^{\theta/2}_T(E_T^{\theta/2}\hat{P}E_T^{\theta/2})^{n-1}E_T^{\theta/2}\tilde{\Pi}(\mathbf{1}\otimes\operatorname{vec}(I_d))}_{\mmu}\\
        =\, & \inner{\mathbf{1}\otimes\operatorname{vec}(I_d)}{\tilde{\Pi}E^{\theta/2}_T(E_T^{\theta/2}\hat{P}E_T^{\theta/2})^{n-1}E_T^{\theta/2}\tilde{\Pi}(\mathbf{1}\otimes\operatorname{vec}(I_d))}_{\mmu},
    \end{aligned}
\end{equation}
where recall $\tilde{\Pi}=\Pi\otimes I_{d^2}$ and $T(x)=\frac{\cos(\phi)}{2}(F(x)\otimes I_{d^2}+I_{d^2}\otimes F(x))$.  The third line holds as $\tilde{\Pi}(\mathbf{1}\otimes \operatorname{vec}(I_d))=\mathbf{1}\otimes \operatorname{vec}(I_d)$ and the last line holds as $\tilde{\Pi}$ is self-adjoint with respect to the inner product.  Next, 
\begin{equation}
    \begin{aligned}
        & \inner{\mathbf{1}\otimes\operatorname{vec}(I_d)}{\tilde{\Pi}E^{\theta/2}_T(E_T^{\theta/2}\hat{P}E_T^{\theta/2})^{n-1}E_T^{\theta/2}\tilde{\Pi}(\mathbf{1}\otimes\operatorname{vec}(I_d))}_{\mmu}\\
        \leq \, & \norm{\tilde{\Pi}E^{\theta/2}_T(E_T^{\theta/2}\hat{P}E_T^{\theta/2})^{n-1}E_T^{\theta/2}\tilde{\Pi}}_{\mmu}\norm{\mathbf{1}\otimes \operatorname{vec}(I_d)}_{\mmu}^2\\
        =\, & d\norm{\tilde{\Pi}E^{\theta/2}_T(E_T^{\theta/2}\hat{P}E_T^{\theta/2})^{n-1}E_T^{\theta/2}\tilde{\Pi}}_{\mmu}\\
        =\, & d\norm{\tilde{\Pi}(\hat{P}E_T^{\theta})^n\tilde{\Pi}}_{\mmu},
    \end{aligned}
\end{equation}
with the last line following from $\tilde{\Pi}\hat{P}=\tilde{\Pi}$ so we can introduce another $\hat{P}$.  Now from Lemma \ref{lem:limit} and the above we have 
\begin{equation}\label{eq:different_norm_bound}
    \begin{aligned}
        \log \norm{E_{T}^{\theta/2}\hat{P}E_{T}^{\theta/2}}_{\mmu} 
        &=\lim_{n\to\infty}\frac{1}{n}\log \inner{\mathbf{1}\otimes\operatorname{vec}(I_d)}{E^{\theta/2}_T(E_T^{\theta/2}\hat{P}E_T^{\theta/2})^{n-1}E_T^{\theta/2}(\mathbf{1}\otimes\operatorname{vec}(I_d))}_{\mmu}\\
        &\leq \lim_{n\to\infty} \frac{1}{n}\log d\norm{\tilde{\Pi}(\hat{P}E_T^{\theta})^n\tilde{\Pi}}_{\mmu}\\
        &=\lim_{n\to\infty} \frac{1}{n}\log \norm{\tilde{\Pi}(\hat{P}E_T^{\theta})^n\tilde{\Pi}}_{\mmu}.
    \end{aligned}
\end{equation}
Although the above may seem a bit contrived, it allows for a general recursive method used in \cite{wagner2006, healy2008, garg2018, qiu2020}; more importantly, we can extend the final result to time-dependent functions via our limit lemma.  In particular, for an initial function $\mathbf{z}_0$, we can track the updates $\mathbf{z}_k:=(\hat{P}E_T^{\theta})^k\mathbf{z}_0$ recursively, and then provide a bound.  This method is strong enough to give the sorts of inequalities we are after.

Now we introduce some notation.  For an element $\mathbf{z}\in\ell_2(\mmu\otimes \mathbf{1})$ let $\mathbf{z}^\parallel=\tilde{\Pi}\mathbf{z}$ and $\mathbf{z}^\perp = (I-\tilde{\Pi})\mathbf{z}$.  In other words, $\mathbf{z}$ decomposes as $\mathbf{z}=\norm{\mathbf{z}^{\parallel}}_{\mmu}(\mathbf{1}\otimes \mathbf{v})+\norm{\mathbf{z}^\perp}_{\mmu}\boldsymbol\rho$ for some $\mathbf{v}$ and some $\boldsymbol\rho$ orthogonal to $\mathbf{1}\otimes \mathbb{C}^{d^2}$.  

As per our discussion above, we will bound the evolution of $\tilde{P}E^{\theta}_T\mathbf{z}$ for any vector $\mathbf{z}$.  We first have the following proposition:
\begin{proposition}
For any $\mathbf{z}\in\ell_2(\mmu\otimes \mathbf{1})$, 
\begin{enumerate}
    \item $\hat{P}\mathbf{z}^\parallel=\mathbf{z}^\parallel$,
    \item $(\hat{P}\mathbf{z})^\perp=\tilde{P}\mathbf{z}^\perp$.
\end{enumerate}
\end{proposition}
\begin{proof}
This follows immediately from the definitions.
\end{proof}

We next have a crucial lemma, the proof deferred to the appendix.  For ease of exposition, we will assume that $\phi=0$ throughout, since the factor of $\cos(\phi)$ is just a scalar and can be pulled through.
\begin{lemma}\label{lem:E_evolve}
For any $\mathbf{z}\in\ell_2(\mmu\otimes \mathbf{1})$, 
\begin{enumerate}
    \item $\norm{(E^{\theta}_T\mathbf{z}^\parallel)^\parallel}_{\mmu}\leq \alpha_1$,
    \item $\norm{(E^{\theta}_T\mathbf{z}^\perp)^\parallel}_{\mmu}\leq \alpha_2$,
    \item $\norm{(E^{\theta}_T\mathbf{z}^\parallel)^\perp}_{\mmu}\leq \alpha_2$,
    \item $\norm{(E^{\theta}_T\mathbf{z}^\perp)^\perp}_{\mmu}\leq \alpha_3$,
\end{enumerate}
where
\[\alpha_1=1+\frac{\mathcal{V}(e^{\mathcal{M} \theta}-\mathcal{M}-\theta-1))}{\mathcal{M}^2}, \alpha_2=\frac{\sqrt{V}(e^{\mathcal{M} \theta}-1)}{\mathcal{M}}, \alpha_3=e^{\mathcal{M} \theta}.\]
\end{lemma}

With this lemma we have the following two claims, with proofs in the appendix.  These are recursive bounds.
\begin{claim}\label{clm:perp_recurse}
    For any $k$, if $\theta<\log(1/\lambda)/\mathcal{M}$, then 
    \[\norm{\mathbf{z}_k^\perp}_{\mmu}\leq \frac{\lambda \alpha_2}{1-\lambda \alpha_3}\max_{0\leq j\leq k}\norm{\mathbf{z}_j^\parallel}_{\mmu}.\]
\end{claim}

\begin{claim}\label{clm:parallel_recurse}
    For any $k$, if $\theta<\log(1/\lambda)/\mathcal{M}$, then 
    \[\norm{\mathbf{z}_k^\parallel}_{\mmu}\leq \bb{\alpha_1+\frac{\lambda \alpha_2^2}{1-\lambda \alpha_3}}\max_{0\leq j\leq k}\norm{\mathbf{z}_j^\parallel}_{\mmu}.\]
\end{claim}

We are now able to prove Theorem \ref{thm:markov_matrix_bernstein}.
\begin{proof}[Proof of Theorem \ref{thm:markov_matrix_bernstein}]
    As in the proof of Theorem \ref{thm:markov_matrix_hoeffding}, Lemma \ref{lem:P_hat_props} and Proposition \ref{prop:T_operator} give
    \begin{equation}
        \mathbb{E}_{\mmu}\bbb{\norm{\prod_{j=1}^n \exp\bb{\frac{\theta e^{i\phi}}{2}F_j(s_j)}}_F^2}
        \leq d\prod_{j=1}^n \norm{E_{T_j}^{\theta/2}\hat{P}E_{T_j}^{\theta/2}}_{\mmu}.
    \end{equation}
    Assume throughout that $\phi=0$ for simplicity, so we can directly use Lemma \ref{lem:E_evolve} without modification.  Equation \ref{eq:different_norm_bound} indicates that 
    \[\log \norm{E_{T_j}^{\theta/2}\hat{P}E_{T_j}^{\theta/2}}_{\mmu}\leq \lim_{n\to\infty}\frac{1}{n}\log \norm{\tilde{\Pi}(\hat{P}E_{T_j}^{\theta})^n \tilde{\Pi}}_{\mmu}.\]
    Claim \ref{clm:parallel_recurse} exactly gives a bound on this term; letting $\mathbf{z}_0\in\ell_2(\mmu\otimes \mathbf{1})$, the claim indicates that 
    \begin{equation}
        \begin{aligned}
            \norm{\tilde{\Pi}(\hat{P}E_{T_j})^n \tilde{\Pi}\mathbf{z}_0}_{\mmu}
            &\leq \bb{\alpha_{j, 1}+\frac{\lambda \alpha_{j, 2}^2}{1-\lambda \alpha_{j, 3}}}^n \norm{\tilde{\Pi}\mathbf{z}_0}_{\mmu}
            &\leq \bb{\alpha_{j, 1}+\frac{\lambda \alpha_{j, 2}^2}{1-\lambda \alpha_{j, 3}}}^n \norm{\mathbf{z}_0}_{\mmu},
        \end{aligned}
    \end{equation}
    with $\alpha_{j, 1}, \alpha_{j, 2}, \alpha_{j, 3}$ as in Lemma \ref{lem:E_evolve} and using $\mathcal{V}_j$ in the definitions of the $\alpha_j$ as we are dealing with time-dependent functions.  Then $\norm{\tilde{\Pi}(\hat{P}E_{T_j}^{\theta})^n \tilde{\Pi}}_{\mmu}\leq \bb{\alpha_{j, 1}+\frac{\lambda \alpha_{j, 2}^2}{1-\lambda \alpha_{j, 3}}}^n$.  Now
    \begin{equation}
    \begin{aligned}
        \bb{\alpha_{j, 1}+\frac{\lambda \alpha_{j, 2}^2}{1-\lambda \alpha_{j, 3}}}^n&=\bb{1+\frac{1}{\mathcal{M}^2}\bb{\mathcal{V}_j(e^{\mathcal{M}\theta}-\mathcal{M}\theta-1)+\frac{\lambda \mathcal{V}_j(e^{\mathcal{M}\theta}-1)^2}{1-\lambda e^{\mathcal{M}\theta}}}}^n \\
        &\leq \exp\bb{\frac{n\mathcal{V}_j}{\mathcal{M}^2}\bb{e^{\mathcal{M}\theta}-\mathcal{M}\theta-1+\frac{\lambda(e^{\mathcal{M}\theta}-1)^2}{1-\lambda e^{\mathcal{M}\theta}}}},
    \end{aligned}
    \end{equation}
    implying 
    \begin{equation}
    \begin{aligned}
        \log \norm{E_{T_j}^{\theta/2}\hat{P}E_{T_j}^{\theta/2}}_{\mmu}
        &\leq \lim_{n\to\infty}\frac{1}{n}\log \norm{\tilde{\Pi}(\hat{P}E_{T_j}^{\theta})^n \tilde{\Pi}}_{\mmu}\\
        &\leq \lim_{n\to\infty}\frac{1}{n}\log \exp\bb{\frac{n\mathcal{V}_j}{\mathcal{M}^2}\bb{e^{\mathcal{M}\theta}-\mathcal{M}\theta-1+\frac{\lambda(e^{\mathcal{M}\theta}-1)^2}{1-\lambda e^{\mathcal{M}\theta}}}}\\
        &=\frac{\mathcal{V}_j}{\mathcal{M}^2}\bb{e^{\mathcal{M}\theta}-\mathcal{M}\theta-1+\frac{\lambda(e^{\mathcal{M}\theta}-1)^2}{1-\lambda e^{\mathcal{M}\theta}}},
    \end{aligned}
\end{equation}
and thus
\begin{equation}
\begin{aligned}
    \mathbb{E}_{\mmu}\bbb{\norm{\prod_{j=1}^n \exp\bb{\frac{\theta }{2}F_j(s_j)}}_F^2}&\leq d\exp\bb{\frac{1}{\mathcal{M}^2}\sum_{j=1}^n \mathcal{V}_j \bb{e^{\mathcal{M}\theta}-\mathcal{M}\theta-1+\frac{\lambda(e^{\mathcal{M}\theta}-1)^2}{1-\lambda e^{\mathcal{M}\theta}}}}\\
    &=d\exp\bb{\frac{\sigma^2}{\mathcal{M}^2} \bb{e^{\mathcal{M}\theta}-\mathcal{M}\theta-1+\frac{\lambda(e^{\mathcal{M}\theta}-1)^2}{1-\lambda e^{\mathcal{M}\theta}}}}
\end{aligned}
\end{equation}
if $\theta<\log(1/\lambda)/\mathcal{M}$, as desired.

Let 
\begin{equation}
M_F(\theta):=d\exp\bb{\frac{\sigma^2}{\mathcal{M}^2} \bb{e^{\mathcal{M}\theta}-\mathcal{M}\theta-1+\frac{\lambda(e^{\mathcal{M}\theta}-1)^2}{1-\lambda e^{\mathcal{M}\theta}}}}.
\end{equation}
To turn this into a tail bound, we will first assume $\mathcal{M}=1$, i.e., $\norm{F_j}\leq 1$ for all $j$, and then drop this assumption later.  Using Equation \ref{eq:concentration_bound_general} we have
\begin{equation}
    \operatorname{Pr}_{\mmu}\bb{\lambda_{\text{max}}\bb{\sum_{j=1}^n F_j(s_j)}\geq t}
    \leq\inf_{\theta>0} d^{1-\pi/4}\exp\bb{-\theta \cdot t +\log M_F\bb{\frac{4}{\pi}\theta}}.
    \end{equation}
Then
\begin{equation}
    \begin{aligned}
        M_F\bb{\frac{4}{\pi}\theta}
        &=d\exp\bb{\sigma^2\bb{e^{4\theta/\pi}-4\theta/\pi-1+\frac{\lambda(e^{4\theta/\pi}-1)^2}{1-\lambda e^{4\theta/\pi}}}}\\
        &=d\exp\bb{\sigma^2\bb{e^{4\theta/\pi}-4\theta/\pi-1+\frac{\lambda(e^{4\theta/\pi}-1)^2}{1-\lambda e^{4\theta/\pi}}}}\\
        &=d\exp\bb{\frac{16}{\pi^2}\bb{\frac{\sigma^2}{16/\pi^2}\bb{e^{4\theta/\pi}-4\theta/\pi-1}+\frac{\sigma^2}{16/\pi^2}\cdot \frac{\lambda(e^{4\theta/\pi}-1)^2}{1-\lambda e^{4\theta/\pi}}}}.
    \end{aligned}
\end{equation}
Write $L=4/\pi$.  Using the inequality $e^x-1\leq xe^x$, the above becomes
\begin{equation}
    \begin{aligned}
        & d\exp\bb{L^2\bb{\frac{\sigma^2}{L^2}\bb{e^{L\theta}-L\theta-1}+\frac{\sigma^2}{L^2}\cdot \frac{\lambda(e^{L\theta}-1)^2}{1-\lambda e^{L\theta}}}}\\
        \leq \, & d\exp\bb{L^2\bb{\frac{\sigma^2}{L^2}\bb{e^{L\theta}-L\theta-1}+\frac{\sigma^2}{L^2}\cdot\frac{\lambda \theta^2e^{2L\theta}}{1-\lambda e^{L\theta}}}}\\
        =\, & d\exp\bb{\sigma^2\bb{e^{L\theta}-L\theta-1}+\frac{\sigma^2\lambda \theta^2e^{2L\theta}}{1-\lambda e^{L\theta}}}.
    \end{aligned}
\end{equation}
Now we argue that 
\[\frac{\sigma^2 \lambda \theta^2 e^{2L\theta}}{1-\lambda e^{L\theta}}\leq \frac{\sigma^2\lambda \theta^2}{1-\lambda-2L\theta}\]
when $0\leq \theta<(1-\lambda)/2L$, which is less than $\log(1/\lambda)$ when $0<\lambda<1$.  Indeed, comparing both sides, it suffices to show
\[\frac{e^{2L\theta}}{1-\lambda e^{L\theta}}\leq \frac{1}{1-\lambda-2L\theta}\]
for this range of $\theta$, which is equivalent to showing $e^{-2L\theta}-\lambda e^{-L\theta}\geq 1-\lambda-2L\theta$.  At $\theta=0$ they are equivalent, so we compare derivatives.  The derivative of the left-hand side is $\lambda Le^{-L\theta}-2Le^{-2L\theta}=Le^{-L\theta}(\lambda-2e^{-L\theta})$, while the derivative of the right-hand side is just $-2L$.  Then
\begin{align*}
    & Le^{-L\theta}(\lambda-2e^{-L\theta})\geq -2L\\
    \iff \, & \frac{e^{-L\theta}(2e^{-L\theta}-\lambda)}{2}\leq 1,
\end{align*}
which is true by seeing that the left-hand side is a decreasing function in $\theta$, so for $\theta\geq 0$ attains its maximum at $\theta=0$.  At this value of $\theta$, the left-hand side is $1-\lambda/2<1$.  

Then
\begin{equation}
    \begin{aligned}
        & d\exp\bb{\sigma^2\bb{e^{L\theta}-L\theta-1}+\frac{\sigma^2\lambda \theta^2e^{2L\theta}}{1-\lambda e^{L\theta}}}\\
        \leq \, & d\exp\bb{\sigma^2\bb{e^{L\theta}-L\theta-1}+\frac{\sigma^2\lambda \theta^2}{1-\lambda-2L\theta}}.
    \end{aligned}
\end{equation}

To recap, at this point we have bounded 
\[M_F(L\theta)\leq d\exp\bb{\sigma^2\bb{e^{L\theta}-L\theta-1}+\frac{\sigma^2\lambda \theta^2}{1-\lambda-2L\theta}}.\]
Then
\begin{align*}
    \operatorname{Pr}_{\mmu}\bb{\lambda_{\text{max}}\bb{\sum_{j=1}^n F_j(s_j)}\geq t}
    & \leq\inf_{\theta>0} d^{1-\pi/4}\exp\bb{-\theta t +\log M_F\bb{\frac{4}{\pi}\theta}}\\
    &\leq d^{2-\pi/4}\inf_{\theta>0}\exp\bb{-\theta t+g(\theta)}
\end{align*}
where 
\[g(\theta):=\sigma^2\bb{e^{L\theta}-L\theta-1}+\frac{\sigma^2 \lambda L^2\theta^2}{1-\lambda-2L\theta}.\]
It is standard that the conjuguate $g^*(t)=\sup_{\theta}\{\theta t-g(\theta)\}$ bounds the tail probability via
\begin{equation}\label{eq:conjugate}
    \operatorname{Pr}\bb{\lambda_{\text{max}}\bb{\sum_{j=1}^n F_j(s_j)}\geq t}\leq d\exp(g^*(t))
\end{equation}
for a convex function $g$ \cite{GORNI1991}.  Using Lemma \ref{lem:conjugate} with $L=4/\pi$, the second statement in Theorem \ref{thm:markov_matrix_bernstein} is immediate in the case $\mathcal{M}=1$.  

More generally, we have
\[
    \operatorname{Pr}\bb{\lambda_{\text{max}}\bb{\sum_{j=1}^n F_j(s_j)}\geq t}
    =\operatorname{Pr}\bb{\lambda_{\text{max}}\bb{\sum_{j=1}^n \frac{F_j(s_j)}{\mathcal{M}}}\geq \frac{t}{\mathcal{M}}}.
\]
Now $\mathbb{E}_{\mmu}[F_j(x)/\mathcal{M}]=0$, $\norm{F_j}\leq 1$, and $\norm{\mathbb{E}_{\mmu}[F_j(x)^2/\mathcal{M}^2]}\leq \mathcal{V}_j/\mathcal{M}^2$.  Plugging this in and rearranging delivers the second statement of the theorem.
\end{proof}

\section{An Application to Covariance Estimation and Markov PCA}\label{sec:pca}
Principal Component Analysis (PCA) is one of the most fundamental problems in data science, used to reduce the dimensionality of multidimensional datasets while retaining as much of the variance as possible and implemented in a wide variety of downstream applications \cite{Jolliffe2002}.  The task is to estimate the largest eigenvector of a population covariance matrix from samples.  Classically, vector-valued samples $\mathbf{x}_1, \dots, \mathbf{x}_n$ are received i.i.d., and the leading eigenvector of the sample covariance matrix is used to estimate the leading eigenvector of the population covariance matrix.

More precisely, we fix a distribution $\mmu$ and let $\Sigma=\mathbb{E}_{\mmu}\bbb{\mathbf{x}\mathbf{x}^\top}\in\bbR^{d\times d}$.  Let $\chi_1\geq  \dots\geq \chi_d$ be the eigenvalues of $\Sigma$ in descending order, and let $\mathbf{v}_1$ be the leading eigenvector of $\Sigma$.  The following is a standard consequence of the classical matrix Bernstein inequality \cite{tropp2015} and Wedin's theorem \cite{wedin1972}:
\begin{theorem}[Offline PCA, \cite{jain2018}]\label{thm:classical_pca}
    Fix $\delta\in (0, 1)$.  Let $\mathbf{x}_1, \dots, \mathbf{x}_n$ be drawn i.i.d. from $\mmu$.  Assume that $\norm{\mathbb{E}_{\mmu}\bbb{\bb{\mathbf{x_j}\mathbf{x_j}^\top -\Sigma}^2}}\leq \mathcal{V}_j$ and $\norm{\mathbf{x}_j\mathbf{x}_j^\top-\Sigma}\leq \mathcal{M}$ almost surely; let $\sigma^2=\sum_{j=1}^n \mathcal{V}_j$.  Let $\hat{\mathbf{v}}$ be the leading eigenvector of $\frac{1}{n}\sum_{j=1}^n \mathbf{x}_j\mathbf{x}_j^\top$.  Then with probability $1-\delta$
    \begin{equation}
        1-\inner{\hat{\mathbf{v}}}{\mathbf{v}_1}^2\leq C_1\frac{\frac{\sigma^2}{n} \ln\bb{\frac{d}{\delta}}}{(\mu_1-\mu_2)^2}\cdot \frac{1}{n}+C_2\bb{\frac{\mathcal{M}\ln\bb{\frac{d}{\delta}}}{\chi_1-\chi_2}}^2\cdot \frac{1}{n^2}
    \end{equation}
    where $C_1, C_2$ are absolute constants.
\end{theorem}

In many data tasks the samples $\mathbf{x}_1, \dots, \mathbf{x}_n$ are not drawn i.i.d. but have instead some dependent structure; for example, this can occur with time series data \cite{hormann2015dynamic, jolliffe2016}.  Specifically, the dependent structure can be Markovian, where the vectors $\mathbf{x}_j:=\mathbf{x}_j(s_j)$ are vector-valued functions of an underlying Markov chain.  An example of this is in the context of token algorithms for Federated PCA, where multiple machines are connected via a graph \cite{grammenos2020federatedprincipalcomponentanalysis, Hartebrodt2021, even2023stochasticgradientdescentmarkovian}: each machine contains some fraction of the total dataset and the goal is to obtain the principal component of the entire dataset with respect to some target stationary distribution; work has also been done to adapt this to the streaming setting \cite{kumar2024} (whereas our following bound is applicable to the offline setting).  The following bound is the adaptation of Theorem \ref{thm:classical_pca} to the Markov setting, which is a consequence of Theorem \ref{thm:markov_matrix_bernstein} and Wedin's theorem; it was first found in \cite{kumar2024} and is derived as an immediate application of our results.  Thus, our theorems are able to directly give the first known bounds for offline Markov PCA in the literature.
\begin{theorem}[Offline Markov PCA, \cite{kumar2024}]\label{thm:markov_pca}
    Fix $\delta \in (0, 1)$.  Let $P$ be a discrete-time Markov chain on general state space with stationary distribution $\mmu$ and absolute spectral gap $\lambda$.  Let $s_1, \dots, s_n$ be a sequence of states driven by $P$ with initial distribution $\mmu$.  Consider a sequence of vector valued functions $\mathbf{x}_j:=\mathbf{x}_j(s_j), j=1, \dots, n$.  Assume that $\norm{\mathbb{E}_{\mmu}\bbb{\bb{\mathbf{x_j}\mathbf{x_j}^\top -\Sigma}^2}}\leq \mathcal{V}_j$ and $\norm{\mathbf{x}_j\mathbf{x}_j^\top-\Sigma}\leq \mathcal{M}$ almost surely; let $\sigma^2=\sum_{j=1}^n \mathcal{V}_j$.  Let $\hat{\mathbf{v}}$ be the leading eigenvector of $\frac{1}{n}\sum_{j=1}^n \mathbf{x}_j\mathbf{x}_j^\top$.  Then with probability $1-\delta$
    \begin{equation}
        1-\inner{\hat{\mathbf{v}}}{\mathbf{v}_1}^2\leq C'_1\frac{\frac{\sigma^2}{n} \ln\bb{\frac{d^{2-\pi/4}}{\delta}}}{(\chi_1-\chi_2)^2}\bb{\frac{1+\lambda}{1-\lambda}}\cdot \frac{1}{n}+C'_2\bb{\frac{\mathcal{M}\ln\bb{\frac{d^{2-\pi/4}}{\delta}}}{(\chi_1-\chi_2)(1-\lambda)}}^2\cdot \frac{1}{n^2}
    \end{equation}
    where $C'_1, C'_2$ are absolute constants.
\end{theorem}
Observe that the bound from Theorem \ref{thm:markov_pca} is a natural generalization of the bound from Theorem \ref{thm:classical_pca} up to constants, the dependence on the absolute spectral gap $\lambda$, and the extra $2-\pi/4$ factor in the dimension (this can be removed with a slightly worse probability of success); when $\lambda=0$, which recovers the independent setting, the bound essentially matches that of Theorem \ref{thm:classical_pca}.  Note that there are now two spectral gaps: the absolute spectral gap $\lambda$ corresponding to the underlying Markov chain, which governs mixing and thus the spectral norm bound between the sample and population covariances, and the difference $\chi_1-\chi_2$ corresponding to the population covariance matrix $\Sigma$ that governs the quality of the estimated eigenvector from Wedin's theorem.  

More generally, our bounds can be useful for analyzing algorithms for covariance matrix estimation for a stationary distribution based on dependent samples from a Markov chain \cite{kook2024covarianceestimationusingmarkov}.

\section{Discussion}
In this work we gave a systematic study of concentration inequalities for sums of Markov-dependent random matrices, namely, sums of time-dependent matrix-valued functions of a nonreversible Markov chain on continuousgi state spaces.  Our techniques broadly follow the classic literature on spectral methods for sums of Markov-dependent random variables, as we bound the largest eigenvalue of certain perturbed operators.  To address time-dependent Markov chains on continuous state spaces we first gave a crucial limit lemma justifying the analysis of the largest eigenvalue of a certain operator, analogous to the scalar case; in the process, we make an interesting connection to the Krein-Rutman theorem, the theory of operators leaving invariant a cone, to tackle the spectral analysis.

In addition to this, we gave a tighter bound on the moment generating function in the Hoeffding setting, exposing the sub-Gaussian nature of the sum, and used this to give improved constants in the final tail bound.  Our technique here was to construct a natural coupling with a two-state chain via a convex majorization argument, similar to the proof of the classical Hoeffding's lemma.  We also gave the first Bernstein-type inequality for sums of Markov-dependent random matrices via a recursive, linear algebraic, argument -- our bound on the moment generating function nicely reveals the Markov dependence, and provides corresponding tail bounds.  Both the Hoeffding and Bernstein inequalities thus generalize exactly, and in the Bernstein case even improves, the scalar setting.  In addition, we expect our general spectral framework to readily extend to other concentration inequalities, such as Bennett's inequality and Bernstein's inequality with unbounded summands (when there is control over the moments).  We hope that our work opens up avenues for application and contributes to a growing body of work on concentration inequalities of sums of random matrices beyond independence.

One interesting direction is a possible improvement of Corollary \ref{cor:nonstationary}.  In general, finite moments of the Radon-Nikodym derivative $\frac{d\boldsymbol{\nu}}{d\mmu}$ can be much smaller than the essential supremum; specifically, it is possible in the scalar setting to tradeoff an improved multiplicative factor involving the $p^{\text{th}}$ moment of the Radon-Nikodym derivative with worse constants in the exponent -- the proof is just H{\"o}lder's inequality.  However, for matrices this argument would force us to consider not the Frobenius norm squared in the moment generating function but the Frobenius norm raised to some power depending on $p$.  It would be interesting if such a tradeoff could be made in our setting as well.

Another interesting direction is with regards to the spectral gap of the chain.  Our ``absolute spectral gap'' is another name for multiplicative reversibilization, commonly used to quantify the mixing of nonreversible chains and is used to reduce analysis to a reversible chain.  Though widespread in the literature, it can be a pessimistic estimate of mixing \cite{chatterjee2023spectral}; it would be nice to have a bound that depends on a more ``nonreversible'' quantity.  In addition, compared to the independent setting, our variance proxy is different and potentially slightly suboptimal; it is in term of a ``sum of norms'' rather than a ``norm of a sum.''  However, we doubt this can be recovered using spectral methods; we believe an improved variance proxy would require new techniques -- such as an application of functional inequalities -- that are not spectral in nature or at least do not rely on the Golden-Thompson inequality utilized here, and thus imply an entirely new approach for ergodic sums of scalar random variables as well.

\section*{Acknowledgments}
JN is supported by NSF CAREER Award 2145800 and NSF DMS-2204449, with partial support by the Deutsche Forschungsgemeinschaft (DFG, German Research Foundation) under Germany’s Excellence Strategy – EXC-2047/1 – 390685813.  BS acknowledges support through NSF IFML grant 2019844.  RW is supported by AFOSR MURI FA9550-19-1-0005, NSF DMS-1952735, NSF DMS-2109155, and NSF IFML grant 2019844.  The authors thank Sanjay Shakkottai for useful references and Purnamrita Sarkar for helpful discussions.


\bibliographystyle{plain}
\bibliography{references}

\renewcommand{\theHsection}{A\arabic{section}}
\begin{appendix}

\section*{\Large Appendix}

\section{Proofs for Section \ref{sec:results} (Main Results)}
\begin{proof}[Proof of Corollary \ref{cor:complex_matrix_concentration}]
    If $Z=X+iY$ is a complex Hermitian $d\times d$ matrix, then $X$ is real symmetric and $Y$ is skew-symmetric.  Then $Z$ is Hermitian if and only if
    \[
    Z'=\begin{bmatrix}
    X & Y\\
    -Y & X
    \end{bmatrix}
    \]
    is a real symmetric $2d\times 2d$ matrix.  If $Z$ has eigenvalues $\lambda_1, \dots, \lambda_d$ then $Z'$ has eigenvalues $\lambda_1, \lambda_1, \dots, \lambda_d, \lambda_d$, i.e., has the same eigenvalues each with twice the multiplicity.  Then call Theorems \ref{thm:markov_matrix_hoeffding} and \ref{thm:markov_matrix_bernstein} with the matrices $Z'$.
\end{proof}
\begin{proof}[Proof of Corollary \ref{cor:nonstationary}]
    Let $H(s_1, \dots, s_n)$ be a measurable function.  Then
    \begin{align*}
        \mathbb{E}_{\boldsymbol\nu}[H(s_1, \dots, s_n)]
        &=\mathbb{E}_{\mmu}\bbb{\frac{d\boldsymbol\nu}{d\mmu}\cdot H(s_1, \dots, s_n)}\\
        &\leq \operatorname{ess\,sup}\frac{d\boldsymbol{\nu}}{d\mmu}\cdot \mathbb{E}_{\mmu}[H(s_1, \dots, s_n)].
    \end{align*}
    This gives the bounds on the moment generating function and the tail bounds then follow.
\end{proof}

\section{Proofs for Section \ref{sec:preliminaries} (Preliminaries)}

\begin{proof}[Proof of Lemma \ref{lem:lambda_Phat}]
The proof follows in a straightforward manner using properties of the operator norm of a tensor product.  We have
\begin{align*}
    \lambda(\tilde{P})
    &=\norm{\tilde{P}-\tilde{\Pi}}_{\mmu}\\
    &=\norm{P\otimes I_{d^2}-\Pi\otimes I_{d^2}}_{\mmu}\\
    &=\norm{(P-\Pi)\otimes I_{d^2}}_{\mmu}\\
    &=\norm{P-\Pi}_{\mmu}\norm{I_{d^2}}\\
    & =\lambda(P).
\end{align*}
\end{proof}

\section{Proofs for Section \ref{sec:operator} (Bounding the operator norm)}
\begin{proof}[Proof of Proposition \ref{prop:T_properties}]
For (1), we have
\begin{align*}
    \mathbb{E}[T(x)]
    &=\frac{\cos(\phi)}{2}(\mathbb{E}[F(x)\otimes I_d]+\mathbb{E}[I_d\otimes F(x)])\\
    &=\frac{\cos(\phi)}{2}(\mathbb{E}[F(x)]\otimes I_d+I_d\otimes \mathbb{E}[F(x)])\\
    &=0,
\end{align*}
since the Kronecker product is a linear operation.

Now we prove (2).  For some $x\in\mathcal{X}$, let $\rho_1(x)\geq \dots\geq \rho_d(x)$ be the eigenvalues of $F(x)$ (the $\rho$ are implicitly scalar-valued functions).  We have that the eigenvalues of $F(x)\otimes I_d$ are exactly these eigenvalues but each with multiplicity $d$ -- this also holds for $I_d\otimes F(x)$.  Now as $F(x)\otimes I_d$ and $I_d\otimes F(x)$ are Hermitian, a crude estimate is that the maximum eigenvalue of the sum is bounded above by $2\rho_1(x)\leq 2b$ and the minimum eigenvalue of the sum is bounded below by $2\rho_d(x)\geq 2a$.  The statement follows.
 
For (3), we have
\[T(x)^2=\frac{\cos^2(\phi)}{4}(F(x)^2\otimes I_d+I_d\otimes F(x)^2+2\cdot F(x)\otimes F(x)).\]
Then
\[\norm{\mathbb{E}[F(x)^2\otimes I_d]}=\norm{\mathbb{E}[F(x)^2]\otimes I_d}\leq \mathcal{V},\]
and the same holds for $\norm{\mathbb{E}[I_d\otimes F(x)^2]}$.  For $F(x)\otimes F(x)$, we consider $(\mathbf{u}\otimes \mathbf{u})^* (F(x)\otimes F(x))(\mathbf{u}\otimes \mathbf{u})$ for some unit vector $\mathbf{u}$.  This is equal to $(\mathbf{u}^*F(x)\mathbf{u})^2$.  Now by Cauchy-Schwartz,
\begin{align*}
    (\mathbf{u}^* F(x)\mathbf{u})^2
    &\leq (\norm{F(x)\mathbf{u}})^2\\
    &=\mathbf{u}^* F(x)^2\mathbf{u}.
\end{align*}
Then $(\mathbf{u}\otimes \mathbf{u})^* \mathbb{E}[F(x)\otimes F(x)](\mathbf{u}\otimes \mathbf{u})\leq \mathbf{u}^* \mathbb{E}[F(x)^2]\mathbf{u}$ for any unit vector $\mathbf{u}$.  It then follows that $\norm{\mathbb{E}[F(x)\otimes F(x)]}\leq \mathcal{V}$.  The statement follows by putting the above together.
\end{proof}

\section{Proofs for Section \ref{sec:final} (Final bounds)}

\subsection{Hoeffding}
\begin{proof}[Proof of Lemma \ref{lem:markov_matrix_convex_majorization}]
    We first note that for any nonnegative convex function $\Psi$, the function $\Psi^p$ is convex for any $p\geq 1$.  We first describe a generating process for $s_j$.  Draw random variables $I_j$ with $I_1=1$ and $I_j\sim  \operatorname{Ber}(1-\lambda)$ for $j>1$ and variables $Z_j\sim \mmu$ uniformly and independent from each other and let
    \begin{equation}
        s_k=\sum_{j=1}^k \bb{\prod_{\ell=j+1}^k (1-I_\ell)}I_j Z_j
    \end{equation}
    so that
    \begin{equation}
        \Psi(F(s_k))=\sum_{j=1}^k \bb{\prod_{\ell=j+1}^k (1-I_\ell)}I_j \Psi(F(Z_j)).
    \end{equation}
    We can verify that this is a valid construction for a Markov chain driven by $\hat{P}$.  We now couple $(y_j, s_j)$ by drawing random variables $B_j\sim \boldsymbol\mu$ independent from each other and let
    \begin{equation}
        y_k=\sum_{j=1}^k \bb{\prod_{\ell=j+1}^k (1-I_\ell)}I_j B_j
    \end{equation}
    so that
    \begin{equation}
        \Psi(G(y_k))=\sum_{j=1}^k \bb{\prod_{\ell=j+1}^k (1-I_\ell)}I_j \Psi(G(B_j)).
    \end{equation}
    We can again verify that this is a valid construction for a Markov chain driven by $\hat{Q}$ and is also a valid coupling.  Then
    \begin{equation}
        \prod_{k=1}^n \Psi(F(s_j))=\prod_{k=1}^n \bb{\sum_{k=1}^n \bb{\prod_{\ell=j+1}^k (1-I_\ell)}I_j \Psi(F(Z_j))}.
    \end{equation}
    We can expand this out as a matrix polynomial.  Each monomial is of degree $n$ and can be determined by a tuple $\alpha\in \{0, 1\}^n$ with $\alpha_1=1$ always, each entry corresponding to $I_j$ and $\Psi(F(Z_j))$.  For each $\alpha$ construct a coefficient $c_\alpha$ and a partition $a$ of $n$ as follows: if $\alpha_j=0$ then $(1-I_j)$ appears in $c_\alpha$ and $a_j=0$.  If $\alpha_j=1$ then $I_j$ appears in $c_\alpha$ and $a_j$ is equal to the one plus the number of zeros appearing sequentially after $\alpha_j$ until the next one.  Then we have $c_\alpha \prod_{j=1}^n \Psi(F(Z_j))^{a_j}$ as our matrix monomial corresponding to $\alpha$.  There are clearly $2^{n-1}$ monomials in this matrix polynomial.

    For example, let $n=6$.  The monomial corresponding to $\alpha=(1, 0, 1, 0, 1, 0)$ is
    \[(1-I_6)I_5(1-I_4)I_3(1-I_2)I_1\Psi(F(Z_5))^2 \Psi(F(Z_3))^2 \Psi(F(Z_1))^2,\]
    the monomial corresponding to $\alpha=(1, 0, 1, 0, 0, 1)$ is 
    \[I_6(1-I_5)(1-I_4)I_3(1-I_2)I_1\Psi(F(Z_6))\Psi(F(Z_3))^3 \Psi(F(Z_1))^2,\]
    the monomial corresponding to $\alpha=(1, 1, 1, 1, 1, 1)$ is
    \[I_6I_5I_4I_3I_2I_1\Psi(F(Z_6))\Psi(F(Z_5))\Psi(F(Z_4))\Psi(F(Z_3))\Psi(F(Z_2))\Psi(F(Z_1)),\]
    and the monomial corresponding to $\alpha=(1, 0, 0, 0, 0, 0)$ is
    \[(1-I_6)(1-I_5)(1-I_4)(1-I_3)(1-I_2)(1-I_1)\Psi(F(Z_1))^6.\]

    Moreover, it is not hard to see that the monomials in this matrix polynomial are pairwise orthogonal (with respect to the trace inner product) using the property that $(1-I_j)I_j=0$.  Therefore
    \begin{equation}
        \begin{aligned}
            \mathbb{E}_{\mmu}\bbb{\norm{\prod_{j=1}^n \Psi(F(s_j))}_F^2}
            &=\mathbb{E}_{\mmu}\bbb{\norm{\sum_\alpha c_\alpha \prod_{j=1}^n \Psi(F(Z_j))^{a_j}}_F^2}\\
            &=\mathbb{E}_{\mmu}\bbb{\sum_\alpha c_\alpha \norm{\prod_{j=1}^n \Psi(F(Z_j))^{a_j}}_F^2}\\
            &=\sum_{\alpha} \mathbb{E}[c_\alpha]\mathbb{E}_{\mmu}\bbb{\norm{\prod_{j=1}^n \Psi(F(Z_j))^{a_j}}_F^2}
        \end{aligned},
    \end{equation}
    where the last line follows from independence of the $I_j$ and $Z_j$.  We now apply Lemma \ref{lem:intermediate_lem} to see that 
    \begin{equation}
        \begin{aligned}
            \mathbb{E}_{\mmu}\bbb{\norm{\prod_{j=1}^n \Psi(F(s_j))}_F^2}
            &=\sum_{\alpha} \mathbb{E}[c_\alpha]\mathbb{E}_{\mmu}\bbb{\norm{\prod_{j=1}^n \Psi(F(Z_j))^{a_j}}_F^2}\\
            &\leq \sum_{\alpha} \mathbb{E}[c_\alpha]\mathbb{E}_{\boldsymbol\mu}\bbb{\norm{\prod_{j=1}^n \Psi(G(B_j))^{a_j}}_F^2}\\
            &=\mathbb{E}_{\boldsymbol\mu}\bbb{\norm{\prod_{j=1}^n \Psi(G(y_j))}_F^2},
        \end{aligned}
    \end{equation} 
    where we use the coupling between $y_j$ and $s_j$.
\end{proof}

\begin{lemma}\label{lem:intermediate_lem}
    For $Z_j$ drawn i.i.d. from $\mmu$ and $B_j$ drawn i.i.d. from $\boldsymbol\mu$, and if $\mathbb{E}_{\mmu}[F]=\mathbb{E}_{\boldsymbol\mu}[G]=0$, and for a nonnegative scalar convex function $\Psi$, we have
    \[\mathbb{E}_{\mmu}\bbb{\norm{\prod_{j=1}^n \Psi(F(Z_j))^{a_j}}_F^2}\leq \mathbb{E}_{\boldsymbol\mu}\bbb{\norm{\prod_{j=1}^n \Psi(G(B_j))^{a_j}}_F^2}.\]
\end{lemma}
\begin{proof}
    We have
    \begin{equation}
    \begin{aligned}
        \mathbb{E}_{\mmu}\bbb{\norm{\prod_{j=1}^n \Psi(F(Z_j))^{a_j}}_F^2}
        &=\mathbb{E}_{\mmu}\bbb{\inner{\prod_{j=1}^n \Psi(F(Z_j))^{a_j}}{\prod_{j=1}^n \Psi(F(Z_j))^{a_j}}}\\
        &=\mathbb{E}_{\mmu}\bbb{\inner{\mathbb{E}_{\mmu}\bbb{\Psi(F(Z_n))^{2a_n}}}{\prod_{j=1}^{n-1}\Psi(F(Z_j))^{a_j}\prod_{j=n-1}^{1}\Psi(F(Z_j))^{a_j}}}
    \end{aligned}
    \end{equation}
using the independence of the $Z_j$.  Now as $\Psi$ is nonnegative and convex, $\Psi^{a_j}$ is convex.  Therefore, by Proposition \ref{prop:matrix_convexity} and the fact that $F$ is mean zero, we have 
\begin{equation}
    \mathbb{E}_{\mmu}\bbb{\Psi(F(Z_n))^{2a_n}}\preceq \mathbb{E}_{\boldsymbol\mu}\bbb{\Psi(G(B_n))^{2a_n}}
\end{equation}
so that
\begin{equation}
    \begin{aligned}
         & \mathbb{E}_{\mmu}\bbb{\inner{\mathbb{E}_{\mmu}\bbb{\Psi(F(Z_n))^{2a_n}}}{\prod_{j=1}^{n-1}\Psi(F(Z_j))^{a_j}\prod_{j=n-1}^{1}\Psi(F(Z_j))^{a_j}}}\\
         \leq \, & \mathbb{E}_{\mmu}\bbb{\inner{\mathbb{E}_{\boldsymbol\mu}\bbb{\Psi(G(B_n))^{2a_n}}}{\prod_{j=1}^{n-1}\Psi(F(Z_j))^{a_j}\prod_{j=n-1}^{1}\Psi(F(Z_j))^{a_j}}}.
    \end{aligned}
\end{equation}
Now since $G$ maps to matrices that are a constant times the identity, the matrices $\mathbb{E}_{\boldsymbol\mu}\bbb{\Psi(G(B_j))^{2a_j}}$ commute with all other matrices.  This allows us to write
\begin{equation}
    \begin{aligned}
        & \mathbb{E}_{\mmu}\bbb{\inner{\mathbb{E}_{\boldsymbol\mu}\bbb{\Psi(G(B_n))^{2a_n}}}{\prod_{j=1}^{n-1}\Psi(F(Z_j))^{a_j}\prod_{j=n-1}^{1}\Psi(F(Z_j))^{a_j}}}\\
        =\, & \mathbb{E}_{\mmu}\bbb{\inner{\mathbb{E}_{\mmu}\bbb{\Psi(F(Z_{n-1}))^{2a_{n-1}}}}{\prod_{j=1}^{n-2}\Psi(F(Z_j))^{a_j}\mathbb{E}_{\boldsymbol\mu}\bbb{\Psi(G(B_n))^{2a_n}}\prod_{j=n-2}^{1}\Psi(F(Z_j))^{a_j}}}.
    \end{aligned}
\end{equation}
We can iterate this process, and using independence of the $B_j$ and rearranging (recall all these resulting matrices commute) we obtain the statement of the lemma.
\end{proof}

\begin{proof}[Proof of Lemma \ref{lem:hoeffding_main}, \cite{fan2021}]
    With some abuse of notation, reuse $a:=\cos(\phi)a$ and $b:=\cos(\phi)b$.  Since $\phi\in[-\pi/2, \pi/2]$, $\phi$ is always nonnegative, so the ordering remains the same.
    
    Let $q=-a/(b-a)$.  The eigenvalues of $K^{\theta}$ are the roots of the quadratic
    \begin{equation}
        \begin{aligned}
            0
            &=\operatorname{det}(\tilde{\eta}I-K^{\theta})\\
            &=\tilde{\eta}^2-(1+\lambda)[(1-p)e^{\theta a}+pe^{\theta b}]\tilde{\eta}+\lambda e^{\theta(a+b)},
        \end{aligned}
    \end{equation}
    where 
    \begin{equation}
        p=\frac{\lambda+(1-\lambda)q}{1+\lambda}.
    \end{equation}
    Then $\eta$ is simply the larger of these two.  Instead of solving directly for $\eta$, we can choose any value that upper bounds $\eta$ that has a considerably simpler form.  In particular, if a quadratic is of the form $x^2-bx+c$, we can choose any $x$ such that $x^2-bx+c\geq 0$ and $x^2\geq c$; it can be verified that such an $x$ upper bounds both roots of the quadratic.  Thus, with our choice of $\tilde{\eta}(\theta)=\exp(\alpha(\lambda)\cdot \theta^2(b-a)^2/8)$, we show
    \begin{equation}\label{eq:quadratic_conditions}
    \begin{aligned}
        & \tilde{\eta}(\theta)^2-(1+\lambda)[(1-p)e^{\theta a}+pe^{\theta b}]\tilde{\eta}(\theta)+\lambda e^{\theta(a+b)}\geq 0,\\
        & \tilde{\eta}(\theta)^2 \geq \lambda e^{\theta(a+b)}.
    \end{aligned}
    \end{equation}
    
    The first point in Equation \ref{eq:quadratic_conditions} holds if and only if 
    \begin{equation}\label{eq:31}
        \frac{\tilde{\eta}(\theta)+\lambda e^{\theta(a+b)}\tilde{\eta}(\theta)^{-1}}{1+\lambda}\geq (1-p)e^{\theta a}+pe^{\theta b}.
    \end{equation}
    Then for the left-hand side of Equation \ref{eq:31} we have
    \begin{align*}
        & \frac{\tilde{\eta}(\theta)+\lambda e^{\theta(a+b)}\tilde{\eta}(\theta)^{-1}}{1+\lambda}\\
        =\, & \frac{\exp(\alpha(\lambda)\theta^2(b-a)^2/8)+\lambda \exp(\theta(a+b)-\alpha(\lambda)\theta^2(a-b)^2/8)}{1+\lambda}\\
        \geq \, & \exp\bb{\frac{\alpha(\lambda)\theta^2(b-a)^2/8+\lambda \theta(a+b)-\lambda \alpha(\lambda)\theta^2(b-a)^2/8}{1+\lambda}}\\
        =\, & \exp\bb{\theta \cdot \frac{\lambda(a+b)}{1+\lambda}+\frac{(b-a)^2\theta^2}{8}\cdot \frac{(1-\lambda)\alpha(\lambda)}{1+\lambda}}\\
        =\, & \exp\bb{\theta \cdot [(1-p)a+pb]+\theta^2\cdot \frac{(b-a)^2}{8}},
    \end{align*}
    where the lower bound is by convexity of the exponential map.  Now the right-hand side of Equation \ref{eq:31} is the moment generating function of the random variable taking value $a$ with probability $1-p$ and value $b$ with probability $p$, and thus is upper bounded by Hoeffding's lemma as
    \[(1-p)e^{\theta a}+pe^{\theta b}\leq \exp\bb{\theta \cdot [(1-p)a+pb]+\theta^2\cdot \frac{(b-a)^2}{8}}.\]
    Putting the two together proves the first line of Equation \ref{eq:quadratic_conditions}.

    The second line of Equation \ref{eq:quadratic_conditions} can be developed as
    \begin{align*}
        \log \bb{\tilde{\eta}(\theta)^2e^{-\theta(a+b)}}
        &=\log \bb{\exp\bb{\frac{\alpha(\lambda)\theta^2(b-a)^2}{4}-\theta(a+b)}}\\
        &=\frac{\alpha(\lambda)\theta^2(b-a)^2}{4}-\theta(a+b).
    \end{align*}
    We show that this is at least $-1/\alpha(\lambda)$.  A sufficient condition for this is
    \[\alpha(\lambda)(a+b)\theta-\frac{\alpha(\lambda)^2 (b-a)^2}{4}\theta^2\leq 1.\]
    The above is maximized at $\theta=2(a+b)/(\alpha(\lambda)(b-a)^2)$, and plugging this in, we have
    \begin{align*}
        \alpha(\lambda)(a+b)\theta-\frac{\alpha(\lambda)^2 (b-a)^2}{4}\theta^2=\frac{(a+b)^2}{(b-a)^2}
    \end{align*}
    which is always at most 1 since $a\leq 0\leq b$.  Then
    \begin{align*}
        \log \bb{\tilde{\eta}(\theta)^2e^{-\theta(a+b)}}
        &=\frac{\alpha(\lambda)\theta^2(b-a)^2}{4}-\theta(a+b)\\
        &\geq -\frac{1}{\alpha(\lambda)}\\
        &=-\frac{1-\lambda}{1+\lambda}\\
        &\geq \log \lambda
    \end{align*}
    when $0\leq \lambda\leq 1$, proving the second line of Equation \ref{eq:quadratic_conditions}.

    As $\tilde{\eta}(\theta)$ has been established to be a suitable upper bound, we then replace $a$ and $b$ in its definition by $\cos(\phi)a$, $\cos(\phi)b$, respectively to obtain the statement of the lemma.
\end{proof}

\subsection{Bernstein}

\begin{proof}[Proof of Lemma \ref{lem:E_evolve}]
    First we use an explicit decomposition of $\mathbf{z}$.  Let $\mathbf{z}=\norm{\mathbf{z}^\parallel}(\mathbf{1}\otimes \mathbf{v}_1)+\norm{\mathbf{z}^\perp}_{\mmu}\boldsymbol\rho$, and $E^\theta_T (\mathbf{1}\otimes \mathbf{v}_1)=a(\mathbf{1}\otimes \mathbf{v}_2)+c\boldsymbol\sigma$, $E^\theta_T \boldsymbol\rho=b(\mathbf{1}\otimes \mathbf{v}_3)+d\boldsymbol\tau$, with $\boldsymbol\rho, \boldsymbol\sigma, \boldsymbol\tau$ norm 1 and orthogonal to $\mathbf{1}\otimes \mathbb{C}^{d^2}$ and $\mathbf{v}_1, \mathbf{v}_2, \mathbf{v}_3$ all norm 1.  Let $D$ be the multiplication operator that acts diagonally by $T$, so that $E^\theta_T=\exp(\theta D)$.
    \begin{enumerate}
        \item We have $\norm{(E^{\theta}_T\mathbf{z}^\parallel)^\parallel}_{\mmu}=|a|\norm{\mathbf{z^\parallel}}_{\mmu}$.  It suffices to bound $|a|$.  First we see that $a=\inner{\mathbf{1}\otimes \mathbf{v}_2}{E^\theta_T (\mathbf{1}\otimes \mathbf{v}_1)}_{\mmu}$.  Then
        \begin{equation}
                \inner{\mathbf{1}\otimes \mathbf{v}_2}{E^\theta_T (\mathbf{1}\otimes \mathbf{v}_1)}_{\mmu}
                =\inner{\mathbf{v}_2}{\mathbb{E}_{\mmu}[\exp(\theta T(x))]\mathbf{v}_1}.
        \end{equation}
          We see this by definition of the multiplication operator $E^\theta_T$ and the definition of the inner product of $\ell_2(\mmu\otimes \mathbf{1})$.  Then
        \begin{equation}
            \begin{aligned}
                \inner{\mathbf{v}_2}{\mathbb{E}_{\mmu}[\exp(\theta T(x))]\mathbf{v}_1}
                &=\mathbb{E}_{\mmu}\bbb{\inner{\mathbf{v}_2}{\exp(\theta T(x))\mathbf{v}_1}}\\
                &=\mathbb{E}_{\mmu}\bbb{\inner{\mathbf{v}_2}{I+\sum_{k=2}^\infty \frac{\theta^k}{k!}T(x)^k \mathbf{v}_1 }}.
            \end{aligned}
        \end{equation}
        We used the fact that $\mathbb{E}_{\mmu}[T(x)]=0$ to remove the $k=1$ term.  Next,
        \begin{equation}
                \mathbb{E}_{\mmu}\bbb{\inner{\mathbf{v}_2}{I+\sum_{k=2}^\infty \frac{\theta^k}{k!}T(x)^k \mathbf{v}_1 }}
                =\inner{\mathbf{v}_2}{\mathbf{v}_1}+\sum_{k=2}^\infty \frac{\theta^k}{k!}\inner{\mathbf{v}_2}{\mathbb{E}_{\mmu}[T(x)^k]\mathbf{v}_1}.
        \end{equation}
        But the inner products in the sum are bounded by the spectral norm of $\mathbb{E}_{\mmu}[T(x)^k]$, so we apply Lemma 6.6.2 from \cite{tropp2015} and Proposition \ref{prop:T_properties}.  Thus
        \begin{equation}
            \begin{aligned}
                \abs{\inner{\mathbf{v}_2}{\mathbb{E}_{\mmu}[\exp(\theta T(x))]\mathbf{v}_1}}
                &\leq 1+\mathcal{V} \sum_{k=2}^\infty \frac{\theta^k}{k!}\mathcal{M}^{k-2}\\
                &=1+\frac{\mathcal{V}(e^{\mathcal{M}\theta}-\mathcal{M}\theta-1)}{\mathcal{M}^2}.
            \end{aligned}
        \end{equation}
        \item We have $\norm{(E^\theta_T \mathbf{z}^\perp)^\parallel}_{\mmu}=|b|\norm{\mathbf{z}^\perp}_{\mmu}$.  It suffices to bound $|b|$.  As before we see $b=\inner{\mathbf{1}\otimes \mathbf{v}_3}{E^\theta_T \boldsymbol\rho}_{\mmu}$.  From our definition of $D$, we have
        \begin{equation}
            \begin{aligned}
                \inner{\mathbf{1}\otimes \mathbf{v}_3}{E^\theta_T \boldsymbol\rho}_{\mmu}
                &=\inner{\mathbf{1}\otimes \mathbf{v}_3}{I+\sum_{k=1}^\infty \frac{\theta^k}{k!}D^k \boldsymbol\rho}_{\mmu}\\
                &=\inner{\mathbf{1}\otimes \mathbf{v}_3}{\boldsymbol\rho}_{\mmu}+\sum_{k=1}^\infty \frac{\theta^k}{k!}\inner{\mathbf{1}\otimes \mathbf{v}_3}{D^k \boldsymbol\rho}_{\mmu}\\
                &=\sum_{k=1}^\infty \frac{\theta^k}{k!}\inner{\mathbf{1}\otimes \mathbf{v}_3}{D^k \boldsymbol\rho}_{\mmu},
            \end{aligned}
        \end{equation}
        since $\mathbf{1}\otimes \mathbf{v}_3$ and $\boldsymbol\rho$ are orthogonal.  Now since $\norm{T}\leq \mathcal{M}$ we also have $\norm{D}_{\mmu}\leq \mathcal{M}$, and so 
        \[-\mathcal{M}^{k-1}\inner{\mathbf{1}\otimes \mathbf{v}_3}{D\boldsymbol\rho}_{\mmu}\leq \inner{\mathbf{1}\otimes \mathbf{v}_3}{D^k \boldsymbol\rho}_{\mmu}\leq \mathcal{M}^{k-1}\inner{\mathbf{1}\otimes \mathbf{v}_3}{D\boldsymbol\rho}_{\mmu}.\]
        This is in turn upper bounded by $\norm{D(\mathbf{1}\otimes \mathbf{v}_3)}_{\mmu}$ and lower bounded by its negation by Cauchy-Schwartz and seeing that $\norm{\boldsymbol\rho}_{\mmu}=1$.  To bound this norm, we have
        \begin{equation}
            \begin{aligned}
                \norm{D(\mathbf{1}\otimes \mathbf{v}_3)}_{\mmu}^2
                &=\inner{D(\mathbf{1}\otimes \mathbf{v}_3)}{D(\mathbf{1}\otimes \mathbf{v}_3)}_{\mmu}\\
                &=\inner{\mathbf{1}\otimes \mathbf{v}_3}{\mathbb{E}_{\mmu}[T(x)^2] (\mathbf{1}\otimes \mathbf{v}_3)}\\
                &\leq \mathcal{V}.
            \end{aligned}
        \end{equation}
        It then follows that
        \begin{equation}
            \begin{aligned}
                \abs{\inner{\mathbf{1}\otimes \mathbf{v}_3}{E^\theta_T \boldsymbol\rho}_{\mmu}}
                &\leq \sqrt{\mathcal{V}}\sum_{k=1}^\infty \frac{\theta^k}{k!}\mathcal{M}^{k-1}\\
                &=\frac{\sqrt{\mathcal{V}}(e^{\mathcal{M} \theta}-1)}{\mathcal{M}}.
            \end{aligned}
        \end{equation}
        \item $\norm{(E^\theta_T \mathbf{z}^\parallel)^\perp}=|c|\norm{\mathbf{z}^\parallel}_{\mmu}$.  Bounding $|c|$, we see $c=\inner{\boldsymbol\sigma}{E^\theta(\mathbf{1}\otimes \mathbf{v}_1)}_{\mmu}$.  Then the result follows as in the previous part.
        \item We simply use the uniform upper bound of $\norm{\exp(\theta T(x))}\leq \exp(\mathcal{M}\theta)$ for all $x\in \mathcal{X}$.  
    \end{enumerate}
\end{proof}

\begin{proof}[Proof of Claim \ref{clm:perp_recurse}]
    We have
    \begin{equation}
        \begin{aligned}
            \norm{\mathbf{z}_k^\perp}_{\mmu}
            &=\norm{\bb{\hat{P}E^\theta_T \mathbf{z}_{k-1}}^\perp}_{\mmu}\\
            &\leq \norm{\bb{\hat{P}E^\theta_T \mathbf{z}_{k-1}^\parallel}^\perp}_{\mmu}+\norm{\bb{\hat{P}E^\theta_T \mathbf{z}_{k-1}^\perp}^\perp}_{\mmu}\\
            &=\norm{\hat{P}\bb{E^{\theta}_T\mathbf{z}_{k-1}^\parallel}^\perp}_{\mmu}+\norm{\hat{P}\bb{E^{\theta}_T\mathbf{z}_{k-1}^\perp}^\perp}_{\mmu}\\
            &\leq \lambda \norm{\bb{E^{\theta}_T\mathbf{z}_{k-1}^\parallel}^\perp}_{\mmu}+\lambda\norm{\bb{E^{\theta}_T\mathbf{z}_{k-1}^\perp}^\perp}_{\mmu}\\
            &\leq \lambda \alpha_2\norm{\mathbf{z}_{k-1}^\parallel}_{\mmu}+\lambda\alpha_3\norm{\mathbf{z}_{k-1}^\perp}_{\mmu}\\
            &\leq \bb{\lambda\alpha_2+(\lambda\alpha_2)(\lambda\alpha_3)+\dots +(\lambda\alpha_2)(\lambda \alpha_3)^{k-1}}\max_{0\leq j\leq k} \norm{\mathbf{z}_j^\parallel}_{\mmu}\\
            &\leq \frac{\lambda \alpha_2}{1-\lambda \alpha_3}\max_{0\leq j\leq k} \norm{\mathbf{z}_j^\parallel}_{\mmu},
        \end{aligned}
    \end{equation}
    where the last step holds if $\lambda \alpha_3=\lambda e^{\mathcal{M}\theta}<1$, implied if $\theta<\log(1/\lambda)/\mathcal{M}$.
\end{proof}

\begin{proof}[Proof of Claim \ref{clm:parallel_recurse}]
    We have
    \begin{equation}
        \begin{aligned}
            \norm{\mathbf{z}_k^\parallel}_{\mmu}
            &=\norm{\bb{\hat{P}E^\theta_T \mathbf{z}_{k-1}}^\parallel}_{\mmu}\\
            &\leq \norm{\bb{\hat{P}E^\theta_T \mathbf{z}_{k-1}^\parallel}^\parallel}_{\mmu}+\norm{\bb{\hat{P}E^\theta_T \mathbf{z}_{k-1}^\perp}^\parallel}_{\mmu}\\
            &\leq \norm{\hat{P}\bb{E^{\theta}_T\mathbf{z}_{k-1}^\parallel}^\parallel}_{\mmu}+\norm{\tilde{P}\bb{E^{\theta}_T\mathbf{z}_{k-1}^\perp}^\parallel}_{\mmu}\\
            &=\norm{\bb{E^{\theta}_T\mathbf{z}_{k-1}^\parallel}^\parallel}_{\mmu}+\norm{\bb{E^{\theta}_T\mathbf{z}_{k-1}^\perp}^\parallel}_{\mmu}\\
            &\leq \alpha_1\norm{\mathbf{z}_{k-1}^\parallel}_{\mmu}+\alpha_2\norm{\mathbf{z}_{k-1}^\perp}_{\mmu}\\
            &\leq \alpha_1\norm{\mathbf{z}_{k-1}^\parallel}_{\mmu}+\frac{\lambda\alpha_2^2}{1-\lambda\alpha_3}\max_{0\leq j\leq k-1}\norm{\mathbf{z}_j^\parallel}_{\mmu}\\
            &\leq \bb{\alpha_1+\frac{\lambda\alpha_2^2}{1-\lambda\alpha_3}}\max_{0\leq j\leq k-1}\norm{\mathbf{z}_j^\parallel}_{\mmu},
        \end{aligned}
    \end{equation}
    where we apply Claim \ref{clm:perp_recurse}.
\end{proof}

The following lemma bounds the conjugate function in Equation \ref{eq:conjugate}.  The proof is adapted from \cite{jiang2018} to our setting and is standard.
    \begin{lemma}\label{lem:conjugate}
    Define 
    \[g_1(t)=\sigma^2\bb{e^{L\theta}-L\theta-1},\quad g_2(t)=\frac{\sigma^2\lambda L^2\theta^2}{1-\lambda-2L\theta}\]
    for any $0\leq \theta<(1-\lambda)/2L$.  If $0<\lambda<1$ then
    \[(g_1+g_2)^*(t)\geq \frac{t^2}{2L^2}\bb{\frac{1+\lambda}{1-\lambda}\cdot \sigma^2+\frac{2Lt}{1-\lambda}}^{-1}.\]
    If $\lambda=0$ then
    \[(g_1+g_2)^*(t)\geq \frac{t^2}{2L^2}\bb{\sigma^2+\frac{Lt}{3}}^{-1}.\]
\end{lemma}
\begin{proof}
    Extend $g_1$ and $g_2$ to all of $\mathbb{R}$:
    \begin{align*}
        & g_1(\theta)=\begin{cases}
        0,\qquad\qquad\qquad\quad\,\, \theta<0\\
        \sigma^2(e^{L\theta}-L\theta-1), \quad  \theta\geq 0
        \end{cases},\\
        & g_2(\theta)=\begin{cases}
            0, \qquad\quad\,\,\,\,\theta<0\\
            \frac{\sigma^2\lambda L^2\theta^2}{1-\lambda-2L\theta},\quad 0\leq \theta<\frac{1-\lambda}{2L},\\
            \infty,\qquad\quad\,\, \theta\geq \frac{1-\lambda}{2L}
        \end{cases}.
    \end{align*}
    These are close and convex and so admit convex conjugates:
    \begin{equation}
        g_1^*(t)=\begin{cases}
            \sigma^2 h_1\bb{\frac{t}{L\sigma^2}},\quad t\geq 0,\\
            \infty, \qquad\qquad\,\,\,\, t<0
        \end{cases},
    \end{equation}
    where $h_1(x)=(1+x)\log(1+x)-x$ and 
    \begin{equation}
        g_2^*(t)=\begin{cases}
            \frac{(1-\lambda)t^2}{2\lambda\sigma^2 L^2}h_2\bb{\frac{2t}{\lambda \sigma^2 L}},\quad t\geq 0,\\
            \infty,\qquad\qquad\qquad\,\,\,\,\, t<0
        \end{cases},
    \end{equation}
    where $h_2(x)=(\sqrt{1+x}+x/2+1)^{-1}$.

    Now $(g_1+g_2)^*$ is consistent with the above definitions, and so $(g_1+g_2)^*(t)=\sup \{\theta t-g_1(\theta)-g_2(\theta)\}$.  If $\lambda>0$, the Moreau-Rockafellar formula \cite{rockafellar1997} states
    \begin{align*}
        (g_1+g_2)^*(t)
        &=\inf\left\{g_1^*(t_1)+g_2^*(t_2)\mid t_1+t_2=t\right\}\\
        &=\inf\left\{\sigma^2 h_1\bb{\frac{t_1}{L\sigma^2}}+\frac{(1-\lambda)t_2^2}{2\lambda\sigma^2 L^2}h_2\bb{\frac{2t_2}{\lambda \sigma^2 L}}\mid t_1+t_2=t, t_1, t_2\geq 0\right\}
    \end{align*}
    Using $h_1(x)\geq x^2/(2(1+x/3))$ and $h_2(x)\geq (2+x)^{-1}$ delivers
    \begin{align*}
        & \inf\left\{\sigma^2 h_1\bb{\frac{t_1}{L\sigma^2}}+\frac{(1-\lambda)t_2^2}{2\lambda\sigma^2 L^2}h_2\bb{\frac{2t_2}{\lambda \sigma^2 L}}\mid t_1+t_2=t, t_1, t_2\geq 0\right\}\\
        \geq \, & \inf\left\{\frac{t_1^2/(2L^2)}{\sigma^2+\frac{Lt_1}{3}}+\frac{t_2^2/(2L^2)}{\frac{2\lambda \sigma^2}{1-\lambda}+\frac{2Lt_2}{1-\lambda}}\mid t_1+t_2=t, t_1, t_2\geq 0\right\}\\
        \geq \, & \inf\left\{\frac{(t_1+t_2)^2/(2L^2)}{\sigma^2+\frac{Lt_1}{3}+\frac{2\lambda \sigma^2}{1-\lambda}+\frac{2Lt_2}{1-\lambda}}\mid t_1+t_2=t, t_1, t_2\geq 0\right\}\\
        =\, & \frac{t^2/(2L^2)}{\frac{1+\lambda}{1-\lambda}\cdot \sigma^2+\frac{2Lt}{1-\lambda}},
    \end{align*}
    where the second to last line uses $t_1^2/a+t_2^2/b\geq (t_1+t_2)^2/(a+b)$ for nonnegative $t_1, t_2$ and positive $a, b$.

    If $\lambda=0$ then $(g_1+g_2)^*(t)=g_1^*(t)$ and follows from the above analogously.
\end{proof}

\end{appendix}

\end{document}